\documentclass[hidelinks,onefignum,onetabnum]{siamart250211}

\usepackage{lipsum}
\usepackage{amsfonts}
\usepackage{graphicx}
\usepackage{epstopdf}
\usepackage[noend]{algorithmic}
\usepackage{caption}
\usepackage{mathtools}
\ifpdf
  \DeclareGraphicsExtensions{.eps,.pdf,.png,.jpg}
\else
  \DeclareGraphicsExtensions{.eps}
\fi

\newsiamremark{remark}{Remark}
\newsiamremark{problem}{Problem}
\newsiamremark{hypothesis}{Hypothesis}
\crefname{hypothesis}{Hypothesis}{Hypotheses}
\newsiamthm{claim}{Claim}

\headers{An \texorpdfstring{$\lowercase{hp}$}{hp} Multigrid Approach for space-time finite elements}{
  Nils Margenberg, Markus Bause, Peter Munch}

\title{An \texorpdfstring{$\lowercase{\boldsymbol{hp}}$}{hp} Multigrid Approach for Tensor-Product Space-Time Finite Element Discretizations of the Stokes Equations\thanks{Submitted to the editors DATE.}}

\author{
  Nils Margenberg\thanks{
    University of Magdeburg,
    Institute for Analysis and Numerics,
    Universit\"atsplatz 2,
    39104 Magdeburg,
    Germany,
    \texttt{nils.margenberg@ovgu.de} (Corresponding
    Author)}
  \and 
  Markus Bause\thanks{Helmut Schmidt University,
    Faculty of Mechanical and Civil Engineering,
    Holstenhofweg 85,
    22043 Hamburg,
    Germany,
    \href{mailto:bause@hsu-hh.de}{\texttt{bause@hsu-hh.de}}}
  \and
  Peter Munch\thanks{Technical University Berlin,
    Faculty of Mathematics,
    Straße des 17. Juni 136,
    10623 Berlin, Germany,
    \href{mailto:p.muench@tu-berlin.de}{\texttt{p.muench@tu-berlin.de}}}}

\usepackage{amsopn}

\renewcommand{\algorithmiccomment}[1]{\bgroup\hfill$\blacktriangleright$~#1\egroup}
\usepackage{diagbox}
\newcommand{\N}{\mathbb{N}} %
\newcommand{\R}{\mathbb{R}} %
\renewcommand{\d}{\mkern3mu\text{d}} %
\usepackage{MnSymbol}
\usepackage[normalem]{ulem}
\usepackage{cancel}
\usepackage{multirow}
\ifpdf
\hypersetup{
  pdftitle={A \texorpdfstring{$hp$}{hp} Multigrid Approach for Tensor-Product Space-Time Finite Element Discretizations of the Stokes Equations},
  pdfauthor={Margenberg, N., Bause, M.,  Munch, P.}
}
\fi

\usepackage{tikz}
\usetikzlibrary{3d,decorations.pathreplacing,shapes, positioning, fit, quotes,angles, fit, calc, patterns}
\usetikzlibrary{spy,shadows}
\tikzset{
  spy using overlaysshadow/.style={
    spy scope={#1,
         every spy on node/.style={
            circle,
            fill, fill opacity=0.25, text opacity=1
             },
         every spy in node/.style={
                 circle, circular drop shadow,
                 fill=white, draw, cap=round
            }
        }
    }
}
\usepackage{pgfplots}
\usepgfplotslibrary{groupplots,colorbrewer,fillbetween}
\pgfplotsset{compat=newest}%
\pgfplotsset{ colormap/Set1-4, cycle multiindex* list={ mark
    list*\nextlist Set1-4\nextlist }, every axis/.append style = {thick},%
}%
\pgfkeys{/pgf/number format/.cd,1000 sep={\,}}
\usepgfplotslibrary{colormaps}
\pgfplotsset{ every mark/.append style={4pt},tick style = {thick,black}}%
\usepackage{booktabs}
\usepackage{enumitem}
\usepackage{siunitx}
\usepackage{fp}
\newcommand{\eval}[2][]{%
    \ifthenelse{\equal{#1}{}}{%
        \FPeval{\result}{#2}%
        \num[round-precision=4]{\result}%
    }{%
        \FPeval{\result}{#2}%
        \num[scientific-notation=false,round-mode=places,print-zero-exponent=false,tight-spacing=true,round-precision=#1]{\result}%
    }%
  }
\usepackage{subcaption}
\definecolor{myblue}{RGB}{0 83 139}
\definecolor{myred}{RGB}{114 16 69}
\definecolor{mygreen}{RGB}{0 94 0}

  \makeatletter
  \newcommand{\customlabel}[2]{%
    \protected@write \@auxout {}{\string \newlabel {#1}{{#2}{\thepage}{#2}{#1}{}} }%
    \hypertarget{#1}{#2}%
  }
  \makeatother

 \newcommand{\MB}[1]{{}{\color{blue} #1}{}}  %
  
\begin{document}

\maketitle
\begin{abstract}
  We present a monolithic \(hp\) space-time multigrid method for tensor-product
  space-time finite element discretizations of the Stokes equations. Geometric
  and polynomial coarsening of the space-time mesh is performed, and the entire
  algorithm is expressed through rigorous mathematical mappings. For the
  discretization, we use inf-sup stable pairs
  \(\mathbb Q_{r+1}/\mathbb P_{r}^{\text{disc}} \) of elements in space and a
  discontinuous Galerkin (DG\( (k) \)) discretization in time with piecewise
  polynomials of order \( k \). The key novelty of this work is the application
  of \( hp \) multigrid techniques in space and time, facilitated and
  accelerated by the matrix-free capabilities of the \texttt{deal.II} library.
  While multigrid methods are well-established for stationary problems, their
  application in space-time formulations encounter unique challenges,
  particularly in constructing suitable smoothers. To overcome these challenges,
  we employ space-time cell and vertex star patch based Vanka smoothers.
  Extensive tests on high-performance computing platforms demonstrate the
  efficiency of our \( hp \) multigrid approach on problem sizes exceeding a
  trillion degrees of freedom (dofs), sustaining throughputs of hundreds of
  millions of dofs per second.
\end{abstract}
\begin{keywords}
  Space-time finite elements, space-time multigrid, monolithic multigrid,
  matrix-free, higher-order finite elements, high-performance computing
\end{keywords}
\begin{MSCcodes}
65M60, 65M55, 65F10, 65Y05
\end{MSCcodes}
\section{Introduction}
Time-dependent partial differential equations, such as the nonstationary Stokes
equations, greatly benefit from methods that exploit parallelism in space and
time. \emph{Space-time finite element methods} (STFEMs) offer a natural
framework for such parallelization by treating time as an additional dimension,
enabling simultaneous discretization and solution in space and time. This work
substantially extends our previous space-time multigrid framework for scalar
parabolic and hyperbolic problems~\cite{margenbergSpaceTimeMultigridMethod2024a}
to the instationary Stokes system. Key advances include support for
vector-valued, multi-variable systems, an integrated $hp$-multigrid hierarchy,
and optional vertex-star patches in the smoother. Together, these developments
yield robust, scalable solvers for large-scale space-time Stokes problems.

In this work, we present a $hp$ space-time multigrid method ($hp$ STMG) for
tensor-product space-time finite element discretizations of the Stokes system.
For the discretization in space we use the mapped version of the inf-sup stable
$\mathbb Q_{r+1}/\mathbb P_{r}^{\text{disc}}$ pairs of finite elements, with
$r\in \N$. For the discretization in time, we use the discontinuous Galerkin
method (DG$(k)$) of order $k\in \N$. Continuous in time Galerkin methods are not
studied here, due to their difficulties related to the computation of an
discrete initial value for the pressure and the non-wellposedness of the
discrete pressure
trajectory~\cite{anselmannOptimalorderPressureApproximation2025}. Other inf-sup
stable element pairs supported by the deal.II library, are also supported by our
implementation. However, we restrict the presentation to the chosen
discretization to maintain focus and avoid unnecessary complexity in the
discussion. While $\mathbb Q_{r+1}/\mathbb P_{r}^{\text{disc}}$ are well-suited
for the cell based Vanka-Smoother, Taylor-Hood elements don't work due to the
larger coupling stencil of the continuous
pressure~\cite{wobkerNumericalStudiesVankaType2009}. Larger patches can
alleviate this problem, see
e.\,g.~\cite{rafiei2025improvedmonolithicmultigridmethods}. The key novelty of
our approach is the application of geometric and polynomial space-time multigrid
techniques. The implementation is facilitated by the matrix-free capabilities of
the \texttt{deal.II}
library~\cite{africa_dealii_2024,kronbichlerGenericInterfaceParallel2012,munchEfficientDistributedMatrixfree2023,fehnHybridMultigridMethods2020}.
The present work builds upon the foundation established
in~\cite{margenbergSpaceTimeMultigridMethod2024a}. The code is available on
GitHub at~\url{https://github.com/nlsmrg/dealii-stfem}.

Extending multigrid methods to STFEMs poses challenges, particularly in the design of effective smoothers for the arising linear systems. However, STFEMs offer appreciable advantages: They naturally integrate spatial and temporal discretizations and handle coupled problems.
Further, they facilitate duality-based and goal-oriented adaptivity in space and
time~\cite{besierGoalorientedSpacetimeAdaptivity2012,bauseFlexibleGoalorientedAdaptivity2021,rothTensorProductSpaceTimeGoalOriented2023}.
Adaptive STFEMs have also been investigated
in~\cite{langerAdaptiveSpaceTime2022,wienersSpacetimeDiscontinuousGalerkin2023,coralloSpaceTimeDiscontinuousGalerkin2023}.
Alternative approaches, particularly for unstructured space-time meshes are
discussed
in~\cite{langerSpaceTimeMethodsApplications2019,steinbachSpaceTimeFiniteElement2015,nochettoSpacetimeMethodsTimedependent2018,ernestiSpaceTimeDiscontinuousPetrov2019a,langerSpaceTimeHexahedralFinite2022}.
The utilization of global STFEMs, i.e.\ the concurrent treatment of all
subintervals, has the potential to fully exploit the computational resources of
high-performance computers. Conversely, local STFEMs employ space-time
variational discretizations as time-marching schemes by selecting a test basis
supported on the subintervals. This requires less computational resources, while
scalability to global formulations is maintained. Finally, higher order FE spaces facilitate improved accuracy of discrete solutions on computationally feasible grids.

Parallel time integration methods have been developed to exploit parallelism in
the time dimension and to overcome the sequential bottleneck of traditional
time-stepping methods. A comprehensive review of such methods can be found
in~\cite{ganderTimeParallelTime2024}. However, most of these methods entail a
trade-off between additional computational complexity and time parallelism.
An alternative is the all-at-once solution of the entire space-time
system~\cite{danieli_space-time_2022,southworth_fast_2022}. Owing to the
established connection between Runge-Kutta methods and variational time
discretizations~\cite{ganderAnalysisNewSpaceTime2016,southworth_fast_2022},
these contributions relate to the present work. Space-time multigrid
methods treat time as an additional grid dimension, enabling simultaneous
multilevel coarsening in space and
time~\cite{hortonSpaceTimeMultigridMethod1995,francoMultigridMethodBased2018,falgoutMultigridMethodsSpace2017b,honBlockToeplitzPreconditioner2023}.
While algebraic multigrid methods have been applied to space-time
systems~\cite{steinbachAlgebraicMultigridMethod2018,langerSpaceTimeHexahedralFinite2022},
geometric multigrid technqiues offer advantages in computational efficiency and
scalability~\cite{hackbuschParabolicMultigridMethods1985,ganderAnalysisNewSpaceTime2016}.
Another approach to time parallelism that does not increase the computational
complexity is stage parallelism within a single time
step~\cite{christliebParallelHighOrderIntegrators2010,paznerStageparallelFullyImplicit2017,munchStageParallelFullyImplicit2023}.
While the scalability of stage parallelism is constrained by the number of
stages, these methods are effective in the scaling limit.

Our goal in this work is the design of a $hp$ STMG with the same
grid-independent convergence seen in established geometric multigrid techniques
for elliptic or stationary Stokes-type
problems~\cite{olshanskiiMultigridAnalysis2012}. For stability reasons, the $hp$
STMG is applied as a preconditioner for GMRES iterations, which has become a
standard approach in multigrid frameworks. For parallel efficiency, the
$V$-cycle form is used, with a single $V$-cycle per application of $hp$ STMG.\@
The efficiency of multigrid methods strongly depend on the smoothing operator.
We employ a space-time cell based Vanka smoother. This choice is motivated by the
proven effectiveness of Additive Schwarz or Vanka-type smoothers in
fluid~\cite{ahmedAssessmentSolversSaddle2018,anselmannGeometricMultigridMethod2023}
and solid mechanics~\cite{wobkerNumericalStudiesVankaType2009}, as well as in
fluid-structure interaction~\cite{failerParallelNewtonMultigrid2021}, dynamic
poroelasticity~\cite{anselmannEnergyefficientGMRESMultigrid2024} and acoustic
wave equations~\cite{margenbergSpaceTimeMultigridMethod2024a}. To ensure
computational efficiency of the $hp$ STMG, we focus on matrix-free
implementations.

In recent years, there has been work on matrix-free, high-order monolithic
multigrid methods to solve Stokes (and Navier-Stokes) equations on
high-performance computing architectures. The authors
of~\cite{kohlTextbookEfficiencyMassively2022} developed a parallel, matrix-free
multigrid solver achieving ``textbook multigrid efficiency'', scaling up to
multiple trillions of unknowns. The authors of~\cite{abu-labdeh_monolithic_2023}
extend a patch-based Vanka smoother to fully implicit Runge-Kutta
discretizations of incompressible flows using standard Taylor-Hood elements.
Jodlbauer et al.\ used a matrix-free monolithic geometric multigrid solver for
discretizations of the Stokes equations with Taylor-Hood elements with scaled
Chebyshev-Jacobi smoothers~\cite{jodlbauerMatrixfreeMonolithicMultigrid2024}.
Prieto Saveedra et al.\ present a matrix-free solver for SUPG and PSPG
stabilized equal-order discretizations of the incompressible Navier-Stokes
equations, and achieve substantial speedups and reduced memory usage compared to
matrix-based methods~\cite{prietosaavedraMatrixFreeStabilizedSolver2024}. The
authors of~\cite{voroninMonolithicMultigridPreconditioners2024} propose a
monolithic $ph$-multigrid method for stationary
Stokes. In line with their findings, our experiments demonstrate that
this approach outperforms geometric multigrid methods.

This paper is organized as follows. In \Cref{sec:stfem} we introduce the
continuous problem and the tensor-product space time finite element
discretization.\@ We formulate the algebraic system arising from the
discretization in \Cref{sec:AlgSys}. In \Cref{sec:mg-framework} we introduce the
$hp$ STMG algorithm, which we use as a preconditioner to a GMRES method. We
verify this methodology by numerical experiments in \Cref{sec:experiments}. We
conclude with an evaluation of the results and a future outlook in
\Cref{sec:conclusions}.

\section{\label{sec:stfem}Continuous and discrete problem}
\subsection{Continuous problem}
We consider the nonstationary Stokes system
\begin{subequations}
  \label{Eq:SE}
  \begin{alignat}{3}
    \label{Eq:SE_1}
    \partial_t \boldsymbol{v} - \nu \boldsymbol\Delta \boldsymbol{v} + \nabla p & = \boldsymbol{f}  && \quad \text{in } \;
    \Omega
    \times (0,\,T)\,,\\
    \label{Eq:SE_2}
    \boldsymbol\nabla \cdot \boldsymbol{v}  & = 0 && \quad \text{in } \; \Omega \times (0,\,T)\,,\\
    \label{Eq:SE_3}
    \boldsymbol{v} (0) & = \boldsymbol{v}_0 && \quad \text{in } \; \Omega \,,\\
    \label{Eq:SE_4}
    \boldsymbol{v}  & = \boldsymbol{0}  && \quad \text{on } \; \partial\Omega \times (0,\,T)\,,
  \end{alignat}
\end{subequations}
where \(\Omega \subset \mathbb{R}^d\), with $d\in \{2,3\}$,  is a bounded open Lipschitz domain and \(T > 0\) is the final time. By $\boldsymbol{v}$ and $p$ we denote the unknown velocity and pressure field, respectively. The force $\boldsymbol{f}$ and initial velocity $\boldsymbol{v}_0$ are prescribed data. In~\eqref{Eq:SE_1}, the coefficient $\nu\in\R_{>0}$ denotes the fluid's viscosity. Homogeneous Dirichlet boundary conditions in~\eqref{Eq:SE_4} are chosen for brevity of presentation. We assume that~\eqref{Eq:SE} admits a sufficiently regular solution up to $t=0$ such that higher order approximations become feasible.

We use standard notation. $H^m(\Omega)$ is the Sobolev space of $L^2(\Omega)$
functions with derivatives up to order $m$ in $L^2(\Omega)$ {while}
$\langle \cdot,\cdot \rangle$ {denotes} the inner product in $L^2(\Omega)$
and its vector-valued and matrix-valued counterparts. Let $L^2_0(\Omega)\coloneq \{ q\in  L^2(\Omega) \mid \int_\Omega q \, \d x =0\}$ and  $H^1_0(\Omega)\coloneq \{u\in H^1(\Omega) \mid u=0 \mbox{ on } \partial \Omega\}$. We put $Q(\Omega)\coloneq L^2_0(\Omega)$ and $\boldsymbol V(\Omega) \coloneq  H^1_0(\Omega)^d$. Here, bold-face letters are used to indicate vector-valued spaces and functions. Further, we define the space
\begin{equation*}
 \boldsymbol V^{\operatorname{div}} (\Omega) \coloneq  \{ \boldsymbol v\in  \boldsymbol V \mid \langle \nabla \cdot  \boldsymbol v, q\rangle = 0 \:\: \forall q\in Q \}\,.
\end{equation*}
For a Banach space $B$ and an interval $J\subset [0,T]$, we let $L^2(J;B)$, $C^m(J;B)$ and $C(J;B)$ denote the Bochner spaces of $B$-valued functions, equipped with their natural norms.  For well-posedness of~\eqref{Eq:SE} in suitable Bochner spaces we refer to, e.\,g.~\cite{john_FiniteElement_2016}.

\subsection{Space-time finite element discretization}
For the discretization of~\eqref{Eq:SE} we use spatial and temporal finite
element meshes, which are combined to a space-time mesh by an algebraic
tensor-product. Discrete space-time function spaces are
then defined in tensor-product form. For the time discretization, we partition
the time interval $I\coloneq(0,\,T]$ into $N$ equal subintervals $I_n
\coloneq(t_{n-1},\,t_n]$, for $n=1,\ldots,N$, where $t_n= n \tau$ and $\tau =
T/N$. Thus, $I=\bigcup_{n=1}^N I_n$. The set $\mathcal{M}_\tau \coloneq \{I_1,\ldots, I_N\}$ of time subintervals is called the time mesh. For $k\in \N_0$, we let $\mathbb P_{k}(J;\R)$ denote the set of all polynomials of degree less than or equal to $k$ on $J\subset I$ with values in $\R$. Then, we put
\begin{equation*}
  \label{Def:Yk}
  Y_\tau^k (\R) \coloneq \left\{w_\tau : I \to \R  \mid w_\tau{}_{|I_n} \in \mathbb
  P_{k}(I_n;\R)\; \forall I_n\in \mathcal M_\tau \right\}\,.
\end{equation*}

For spatial discretization, let $\mathcal{T}_h$ be a shape-regular
triangulation of $\Omega$ into quadrilateral and hexahedral elements in two and three space
dimensions with mesh size
$h>0$. These element types are chosen for our implementation that is based on
the deal.II library~\cite{africa_dealii_2024}.
We define the local finite element spaces on each cell \( K \in \mathcal{T}_h \)
by
\begin{equation}\label{Def:VKQK}
  \boldsymbol{V}_{r+1}(K) \coloneq \left( \mathbb{Q}_{r+1} \right)^d \circ \boldsymbol T_K, \quad Q_r(K) \coloneq \mathbb{P}_r^{\text{disc}} \circ \boldsymbol T_K,
\end{equation}
for some $r\geq 1$; cf.~\cite[Section 3.64]{john_FiniteElement_2016}. Here
\( \boldsymbol T_K \) denotes the standard multilinear mapping of polynomials on
the reference element. We note that either the mapped or unmapped variant of
\( \mathbb{P}_r^{\text{disc}} \) may be used. We note that either the mapped or
unmapped variant of \( \mathbb{P}_r^{\text{disc}} \) may be used. We choose the
mapped variant, as it ensures geometric consistency on curved or non-affine
meshes, and leads to better-conditioned system matrices. The unmapped variant,
while simpler to evaluate, may deteriorate the conditioning of the
discretization.
Based on~\eqref{Def:VKQK}, the global finite element
spaces for approximating $\boldsymbol V(\Omega)$ and $Q(\Omega)$ and defining the $hp$ STMG are
\begin{subequations}
  \label{Def:VhQh}
  \begin{alignat}{2}
    \label{Def:Vh}
    \boldsymbol  V_h^{r+1} (\Omega) & \coloneq \{\boldsymbol v_h \in \boldsymbol V \; : \; \boldsymbol v_{h}{}_{|K}\in  {\boldsymbol V_{r+1}(K)} \;\; \text{for all}\; K \in \mathcal{T}_h\}\,, \\[3pt]
    \label{Def:Qh}
    Q_h^r (\Omega) & \coloneq \{q_h \in Q \; : \; q_{h}{}_{|K}\in {Q_r(K)} \;\; \text{for all}\; K \in \mathcal{T}_h\}\,,\\[3pt]
    \label{Def:Qhp}
    Q_h^{r,+}(\Omega) & \coloneq Q_h^r (\Omega) \oplus \operatorname{span}\{1\}\,.
  \end{alignat}
\end{subequations}
The definition of $Q_h^r$ leads to a discontinuous (in space) pressure
approximation. Further, $Q_h^{r,+}(\Omega)$ is the pressure finite element space
without orthogonality condition. These spaces are used to define the grid transfer process of the multigrid approach and ensure that the discrete pressure is kept in the correct space; cf.~\eqref{Eq:DefRsSs}. The space of discretely divergence-free functions is given by
\begin{equation}
  \label{Def:Vdiv}
  \boldsymbol V_h^{\operatorname{div}} (\Omega)\coloneq \{\boldsymbol v_h \in \boldsymbol V_h^{r+1}(\Omega) \mid \langle \nabla \cdot \boldsymbol v_h,q_h\rangle  = 0 \; \text{for all } q_h \in Q_h^r(\Omega)\}\,.
\end{equation}
\begin{remark}[Choice of the finite element spaces in~\eqref{Def:VhQh}]
  \begin{itemize}[leftmargin=*,itemsep=0pt,topsep=0pt,parsep=0pt]
  \item The \( \mathbb{Q}_{r+1}/\mathbb{P}_r^{\text{disc}} \) pair is
    particularly well suited for the cell based Vanka smoother due to its local
    coupling structure. While this has been confirmed empirically in several
    publications
    (cf.~\cite{ahmedAssessmentSolversSaddle2018,john_numerical_2000,turek_efficient_1999,wesseling_geometric_2001}
    and references therein), a rigorous theoretical analysis appears to be
    lacking in the literature. Other inf-sup stable element pairs, such as
    Taylor-Hood elements, are supported by our implementation. Due to the
    continuous pressure space, Taylor-Hood elements induce a larger coupling
    stencil and require overlapping patches (e.\,g. vertex star patches) for
    effective smoothing.
  \item Although \( \mathbb{P}_r^{\text{disc}} \) elements are not
    tensor-product spaces, truncated tensor-product bases in \texttt{deal.II}
    allow them to be used efficiently within the matrix-free framework
    (cf.~\Cref{sec:mf}).
  \end{itemize}
\end{remark}

The global fully discrete solution spaces are now defined by the tensor-products
\begin{equation}
\label{Eq:GDS}
\boldsymbol  H_{\tau,h}^{\boldsymbol v} =   Y_\tau^k (I) \otimes \boldsymbol  V_h^{r+1}(\Omega) \,, \quad
H_{\tau,h}^{p}  =   Y_\tau^k (I) \otimes \boldsymbol  Q_h^r(\Omega)\,.
\end{equation}
\begin{remark}[Function spaces and their tensor product structure]
\begin{itemize}[leftmargin=*,itemsep=0pt,topsep=0pt,parsep=0pt]
\item The algebraic tensor product $Y_\tau(I) \otimes V_h(\Omega)$ of two finite element spaces $Y_\tau(I)$ and $V_h(\Omega)$ is defined by
\begin{equation*}
Y_\tau(I) \otimes V_h(\Omega) \coloneq \operatorname{span}\{f\otimes g \mid f \in Y_\tau(I)\,, \; g \in V_h(\Omega)\}	\,,
\end{equation*}
with mapping $f\otimes g: (t,\boldsymbol x) \to f(t)g(\boldsymbol x)$. For the construction principle of tensor products of Hilbert spaces we refer to~\cite[Section 1.2.3]{picard_partial_2011}.
~
\item The Hilbert spaces $\boldsymbol  H_{\tau,h}^{\boldsymbol v}$ and $H_{\tau,h}^{p}$ are isometric to the Bochner spaces $Y_\tau^k(I;$ $\boldsymbol V_h^{r+1}(\Omega))\subset L^2(I;\boldsymbol V^{r+1}(\Omega))$ and $Y_\tau^k(I;Q_h^{r}(\Omega))\subset L^2(I;Q)$ of piecewise polynomials with values in $\boldsymbol V_h^{r+1}(\Omega)$ and $Q_h^{r}(\Omega)$; cf.~\cite[Proposition~1.2.28]{picard_partial_2011}.
\end{itemize}
\end{remark}
For any function $w: I\to \boldsymbol V_h$ that is piecewise sufficiently smooth with respect to the time mesh $\mathcal{M}_{\tau}$, for instance for $w\in \boldsymbol  H_{\tau,h}^{\boldsymbol v} $, we define the right-hand side limit at a point $t_n$ by $w^+(t_n) \coloneq \lim_{t\to t_n+0} w(t)$ for $0\leq n<N$. Now, we introduce the fully discrete space-time finite element approximation of~\eqref{Eq:SE}.

\begin{problem}[Discrete variational problem]\label{Prob:DVP}
Let the data $\boldsymbol f\in L^2(I;L^2(\Omega))$ and an approximation $\boldsymbol{v}_{0,h}\in V_h^{\operatorname{div}} (\Omega)$ of $\boldsymbol{v}_0\in \boldsymbol{V}^{\operatorname{div}} (\Omega) $ be given. Put  $\boldsymbol v_{\tau,h}(t_0)\coloneq \boldsymbol{v}_{0,h}$. Find $(\boldsymbol v_{\tau,h},p_{\tau,h})\in \boldsymbol  H_{\tau,h}^{\boldsymbol v} \times H_{\tau,h}^{p}$ such that for all $(\boldsymbol{w}_{\tau,h},q_{\tau,h}) \in \boldsymbol  H_{\tau,h}^{\boldsymbol v} \times H_{\tau,h}^{p}$ there holds that
\begin{subequations}
\label{Eq:DSE_0}
\begin{alignat}{2}
\nonumber
\sum_{n=1}^N \int_{t_{n-1}}^{t_n} \langle \partial_t \boldsymbol v_{\tau,h}, \boldsymbol w_{\tau,h} \rangle +  \nu \langle \nabla \boldsymbol v_{\tau,h}, & \nabla \boldsymbol w_{\tau,h}\rangle -  \langle p_{\tau,h}, \nabla \cdot \boldsymbol w_{\tau,h} \rangle \d t   \\
\label{Eq:DSE_1}
  + \sum_{n=0}^{n-1} \langle \lsem \boldsymbol v_{\tau,h} \rsem_n, \boldsymbol w_{\tau,h}^+(t_n)\rangle & =  \sum_{n=1}^N \int_{t_{n-1}}^{t_n} \langle \boldsymbol f, \boldsymbol w_{\tau,h} \rangle \d t\,, \\[3pt]
\label{Eq:DSE_2}
\sum_{n=1}^N \int_{t_{n-1}}^{t_n} \langle \nabla \cdot \boldsymbol v_{\tau,h}, q_{\tau,h} \rangle \d t& = 0\,,
\end{alignat}
\end{subequations}
with the jump $\lsem \boldsymbol v_{\tau,h}\rsem_n \coloneq  \boldsymbol v_{\tau,h}^+(t_n) - \boldsymbol v_{\tau,h}(t_n)$.
\end{problem}
Well-posedness of Problem~\ref{Prob:DVP} can be shown along the lines of~\cite[Lemma 3.2]{anselmannEnergyefficientGMRESMultigrid2024}.

\section{\label{sec:AlgSys}Algebraic system}
Here, we rewrite Problem~\ref{Prob:DVP} in its algebraic form by exploiting the
tensor product structure~\eqref{Eq:GDS} of the discrete spaces. In \Cref{sec:mg-framework} we then embed the algebraic system into an $hp$ multigrid approach.

\subsection{Preliminaries}
\label{Sec:Prem}

For time integration in~\eqref{Eq:DSE_0}, it is natural to apply the right-sided $(k+1)$-point Gau{ss}--Radau quadrature formula. On $I_n$, it reads as
\begin{equation}
  \label{Eq:GF}
  Q_n(w) \coloneq \frac{\tau_n}{2}\sum_{\mu=1}^{k+1} \hat
  \omega_\mu w(t_n^{\mu}) \approx \int_{I_n} w(t) \d t \,,
\end{equation}
where $t_n^{\mu} =T_n(\hat t_{\mu})$, for $\mu = 1,\ldots,k+1$, are the
Gauss--Radau quadrature  points on $I_n$ and $\hat \omega_\mu$ the corresponding
weights. Here, $T_n(\hat t)\coloneq(t_{n-1}+t_n)/2 + (\tau_n/2)\hat t$ is the affine
transformation from $\hat I = [-1,1]$ to $I_n$ and $\hat t_{\mu}$ are the
Gau{ss}--Radau quadrature points on $\hat I$. The quadrature rule~\eqref{Eq:GF} is exact for all $w\in \mathbb P_{2k} (I_n;\R)$, and $t_n^{k+1}=t_n$.

For time interpolation, a Lagrangian basis with respect to the Gauss--Radau quadrature points and with local support on the subintervals $I_n$, for $n=1,\ldots, N$, is used,
\begin{equation}
\label{Eq:BasYk}
\begin{aligned}
Y_\tau^k(I) = \operatorname{span}\big\{& \varphi_{n}^a \in L^2(I) \mid \varphi^a_{n}{}_{|I_b}\in \mathbb P_k(I_b;\R)\,, \text{ for } b=1,\ldots,N\,, \; \\
&  \operatorname{supp}\, \varphi^a_{n}\subset \overline I_n\,,\;  \varphi^a_{n} (t_{n}^{\mu}) = \delta_{a,\mu} \,,\text{ for } \mu =1,\ldots, k+1\,, \\
&   \text{ and for } a=1,\ldots,k+1\,,  \; n=1,\ldots, N \big \}\,,
\end{aligned}
\end{equation}
with the Kronecker symbol $\delta_{a,\mu}$. For space discretization, we use the standard (global) finite element bases associated with the spaces~\eqref{Def:VhQh} and put
\begin{subequations}
\label{Eq:BasVhQh}
\begin{alignat}{2}
\label{Eq:BasVh}
\boldsymbol  V_h^{r+1} (\Omega) & =  \operatorname{span}\big\{\boldsymbol \chi_m^{\boldsymbol v} \mid m=1,\ldots ,M^{\boldsymbol v}  \big\}\\[3pt]
\label{Eq:BasQh}
Q_h^{r,+} (\Omega) & =  \operatorname{span}\big\{\chi^{p}_m \mid m=1,\ldots ,M^p  \big\}\,.
\end{alignat}
\end{subequations}
Then, functions $(\boldsymbol v_{\tau,h},p_{\tau,h}) \in \boldsymbol  H_{\tau,h}^{\boldsymbol v} \times H_{\tau,h}^{p}$ admit for $\boldsymbol{x}\in \Omega$ and $t\in I$ the representation
\begin{subequations}
\label{Eq:Repvp}
\begin{alignat}{2}
\label{Eq:Repv}
\boldsymbol v_{\tau,h}(\boldsymbol x,\,t) &  = \sum_{n=1}^N \sum_{a=1}^{k+1} \sum_{m=1}^{M^{\boldsymbol v}}  v_n^{a,m} \varphi^a_{n}(t) \boldsymbol  \chi_m^{\boldsymbol v}(\boldsymbol x)\,,\\[3pt]
\label{Eq:Repp}
p_{\tau,h}(\boldsymbol x,\,t) & =  \sum_{n=1}^N \sum_{a=1}^{k+1} \sum_{m=1}^{M^{p}} p_n^{a,m}\varphi^a_{n}(t) \chi_m^{p}(\boldsymbol x)
\end{alignat}
\end{subequations}
with coefficients $v_n^{a,m} \in \R^d$ and $p_n^{a,m}\in \R$ for $n=1,\ldots,N$,
$a=1,\ldots,k+1$ and $m=1,\ldots, M$, with $M\in \{M^{\boldsymbol v},M^p\}$.
This representation again shows the tensor-product structure, which we exploit
here. We note that the orthogonality condition for the pressure is not implemented yet in~\eqref{Eq:Repp}, this will be done below in the multigrid framework in \Cref{Sec:GTO}.

To recast~\eqref{Eq:DSE_0} in algebraic form, we use a local (i.e. on each subinterval $I_n$) space and variable major order of the coeffients $v_n^{a,m} \in \R$ and $p_n^{a,m}\in \R$. For this, we introduce the column vectors
\begin{equation}
\label{Eq:DefSubvec_1}
\boldsymbol V_{n}^{a}  \coloneq \big(v_{n}^{a,1},\ldots,v_{n}^{a,{M^{\boldsymbol v}} }\big)^\top\in \R^{{M^{\boldsymbol v}} } \,,\quad\boldsymbol P_{n}^{a}   \coloneq \big(p_{n}^{a,1},\ldots,p_{n}^{a,M^p}\big)^\top\in \R^{M^p}\,,
\end{equation}
for $a=1,\ldots,k+1$. From~\eqref{Eq:DefSubvec_1} we define the column vectors
\begin{equation}
\label{Eq:DefSubvec_2}
  \boldsymbol V_{n}  \coloneq\big(\boldsymbol V_{n}^{1},\ldots ,\boldsymbol V_{n}^{k+1}\big)^\top\in \R^{(k+1)\cdot {M^{\boldsymbol v}} } \,,\quad
  \boldsymbol  P_{n}  \coloneq \big(\boldsymbol  P_{n}^{1},\ldots,\boldsymbol  P_{n}^{k+1}\big)^\top\in \R^{(k+1) \cdot M^p}
\end{equation}
for $n=1,\ldots,N$. For improved readability, the transpose sign is
skipped for the subvectors $\boldsymbol V_{n}^{a}$ and $\boldsymbol P_{n}^{a}$
in~\eqref{Eq:DefSubvec_2}. Throughout the paper, we don't differ in the notation between column and row vectors, if the meaning is clear from the context.
The global column vector $\boldsymbol  X $ of unknowns on $\Omega\times I$, with $\boldsymbol  X_n=(\boldsymbol  V_n,\boldsymbol  P_n)^\top\in \R^{(k+1)\cdot ({M^{\boldsymbol v}}  + M^p)}$ for $n=1,\ldots, N$,  is then defined by
\begin{equation}
  \label{Eq:DefX}
  \boldsymbol  X = (\boldsymbol  X_1,\ldots,\boldsymbol  X_N)^\top \coloneq \big(\boldsymbol  V_{1},\boldsymbol P_{1},\ldots, \boldsymbol  V_{N},\boldsymbol P_{N}\big)^\top\in \R^{N\cdot (k+1)\cdot ({M^{\boldsymbol v}}  + M^p)}\,.
\end{equation}
For the temporal finite element basis induced by~\eqref{Eq:BasYk}, we define the local matrices $\boldsymbol K^\tau_n\in \R^{(k+1),(k+1)}$, $\boldsymbol M_n^\tau\in \R^{(k+1),(k+1)}$ and $\boldsymbol C^\tau_n\in \R^{(k+1),(k+1)}$  by
\begin{subequations}
\label{Eq:DefKMC}
\begin{alignat}{2}
  (\boldsymbol K^\tau_n)_{a,b} & \coloneq \int_{t_{n-1}}^{t_n} \partial_t \varphi_n^b(t) \, \varphi_n^a(t)\d t + \varphi_n^b (t_{n-1}^+) \,  \varphi_n^a (t_{n-1}^+) \,, \\[3pt]
  (\boldsymbol M_n^\tau)_{a,b} &   \coloneq \int_{t_{n-1}}^{t_n} \varphi_n^b(t) \, \varphi_n^a(t) \d t\,, \\[3pt]
  (\boldsymbol C^\tau_n)_{a,b} &  \coloneq
  \left\{\begin{array}{@{}ll}
  \varphi_{n-1}^b (t_{n-1}) \, \varphi_n^a (t_{n-1}^+)\,, & \text{for } n>1 \,, \\[3pt]
  \left.\begin{array}{@{}lll@{}}
    \varphi_n^a (t_{n-1}^+)\,, & \text{for } b= k+1\,,\\[3pt]
  0 \,, & \text{for } b\in \{1,\ldots, k\} \,,
  \end{array}\right\} & \text{for } $n=1$\,,
  \end{array}\right.
\end{alignat}
\end{subequations}
for $a,b=1,\ldots,k+1$. For the spatial finite element basis induced by~\eqref{Eq:BasVhQh},  we let $\boldsymbol  M_h^{\boldsymbol  v} \in \R^{M^{\boldsymbol  v},M^{\boldsymbol  v}}$, $\boldsymbol A_h \in \R^{M^{\boldsymbol  v}, M_h^{\boldsymbol  v}}$, $\boldsymbol  B \in \R^{M^p, M^{\boldsymbol  v}}$ and $\boldsymbol  M_h^{p} \in \R^{M^{p},M^{p}}$ be defined by
\begin{subequations}
\label{Eq:DefMAB}
\begin{alignat}{4}
\label{Eq:DefMAB_1}
(\boldsymbol M_h)_{i,j} &  \coloneq  \int_\Omega \boldsymbol  \chi_j^{\boldsymbol v}(\boldsymbol x) \,  \boldsymbol  \chi_i^{\boldsymbol v}(\boldsymbol x)  \d \boldsymbol x\,, &
(\boldsymbol A_h)_{i,j} &  \coloneq  \int_\Omega \nabla \boldsymbol  \chi_j^{\boldsymbol v}(\boldsymbol x) \cdot \nabla \boldsymbol  \chi_i^{\boldsymbol v}(\boldsymbol x)  \d \boldsymbol x\,, \\[3pt]
\label{Eq:DefMAB_2}
(\boldsymbol B_h)_{l,j}  &  \coloneq  \int_\Omega \nabla \cdot \boldsymbol  \chi_j^{\boldsymbol v}(\boldsymbol x) \,  \chi_i^{p}(\boldsymbol x)   \d \boldsymbol{x}\,,\quad &
(\boldsymbol M^p_h)_{l,m} &  \coloneq  \int_\Omega \chi_m^{p}(\boldsymbol x) \,  \chi_l^{p}(\boldsymbol x)  \d \boldsymbol x
\end{alignat}
\end{subequations}
for $i,j=1,\ldots,M^{\boldsymbol{v}}$ and $l,m=1,\ldots,M^p$.
Next, we introduce the right-hand side column vector
\begin{equation}
  \label{Eq:DefB}
  \boldsymbol  B = \big(\boldsymbol  B_1, \ldots,  \boldsymbol  B_N\big)^\top\in \R^{N\cdot (k+1)\cdot ({M^{\boldsymbol v}}  + M^p)}\,,
  \quad \text{with}\; \;
  \boldsymbol  B_n = (\boldsymbol F_n,\boldsymbol{0})^\top
\end{equation}
for $n=1,\ldots,N$ and subvectors $\boldsymbol F_{n}$ defined by
\begin{equation}
\label{Eq:DefFn}
\boldsymbol F_{n} \coloneq (\boldsymbol F_{n}^{1},\ldots,\boldsymbol F_{n}^{k+1})^\top \in \R^{(k+1)\times M^{\boldsymbol v}}\,, \quad \text{with}\; \;
(\boldsymbol F_{n,}^{a})_i \coloneq Q_n(\langle \boldsymbol f, \varphi_n^a\, \boldsymbol{\chi}_i^{\boldsymbol{v}} \rangle)
\end{equation}
for $a=1,\ldots,k+1$ and $i=1,\ldots,M^{\boldsymbol{v}}$, and with the quadrature formula~\eqref{Eq:GF}. For the well-definedness of $(\boldsymbol F_{n,}^{a})_i $ in~\eqref{Eq:DefFn} we tacitly make the stronger regularity assumption that $\boldsymbol{f}\in C(I;\boldsymbol{L}^2(\Omega))$ is satisfied.

Finally, we recall the tensor (or right Kronecker) product  $\boldsymbol A\otimes \boldsymbol B$ of matrices $\boldsymbol A\in \R^{r,r}$ and $\boldsymbol B\in \R^{s,s}$, for $r,s\in \N$, defined by
\begin{equation}
  \label{Eq:DefKr}
  \boldsymbol A \otimes \boldsymbol B \coloneq \begin{pmatrix}
    a_{1,1} \boldsymbol B & \cdots & a_{1,r}\boldsymbol B\\
    \vdots  & \ddots & \vdots\\
    a_{r,1}\boldsymbol B & \cdots & a_{r,r}\boldsymbol B
  \end{pmatrix} = \left(a_{ij}\boldsymbol B\right)_{i,j=1}^{r}\,.
\end{equation}

\subsection{Algebraic form of the discrete problem}
\label{Sec:Afdp}

In Problem~\ref{Prob:DVP} we choose a tensor product basis of the solution and test space, with the natural Lagrangian basis of~\eqref{Eq:BasYk} built of functions supported on a single subinterval $I_n$. Then we recast for~\eqref{Eq:DSE_0} the following sequence of local problems on $I_n$.
\begin{problem}[Local algebraic problem]
\label{Prob:LocAlg}
Let $n\in \{1,\ldots,N\}$. For $n>1$ let $\boldsymbol{v}_{\tau,h}(t_{n-1})=\sum_{m=1}^{M^{\boldsymbol v}} v_{n-1}^{k+1,m} \boldsymbol{\chi}_m^{\boldsymbol v}$. For $n=1$ and $\boldsymbol v_{0,h}\in \boldsymbol{V}_h^{\operatorname{div}}$ let $\boldsymbol{v}_{0,h}=\sum_{m=1}^{M^{\boldsymbol v}} v_0^{m} \, \boldsymbol{\chi}_m^{\boldsymbol v}$. Put
\begin{equation}
\label{Eq:DefVnm1}
\boldsymbol{V}_{n-1}\coloneq
\left\{\begin{array}{@{}ll}
\big(\boldsymbol 0, \ldots, \boldsymbol 0,v_{n-1}^{k+1,1},\ldots,v_{n-1}^{k+1,M^{\boldsymbol v}} \big)^\top\,, & \text{for } n>1\,,\\[3pt]
\big(\boldsymbol 0, \ldots, \boldsymbol 0,v_{0}^{1},\ldots,v_{0}^{M^{\boldsymbol v}} \big)^\top\,, & \text{for } n=1\,.
\end{array} \right.
\end{equation}
Find $(\boldsymbol{V}_n,\,\boldsymbol{P}_n)\in \R^{(k+1)(M^{\boldsymbol{v}}+M^p)}$ such that
\begin{equation}
\label{Eq:InAlg}
\begin{pmatrix}
\boldsymbol K_n^\tau \otimes \boldsymbol M_h + \boldsymbol M_n^\tau \otimes \boldsymbol A_h  &  \boldsymbol M_n^\tau \otimes \boldsymbol B_h^\top \\[3pt]
\boldsymbol M_n^\tau \otimes \boldsymbol B_h & \boldsymbol 0
\end{pmatrix}
\begin{pmatrix}
\boldsymbol{V}_n\\[3pt] \boldsymbol{P}_n
\end{pmatrix}
=
\begin{pmatrix}
\boldsymbol{F}_n\\[3pt] \boldsymbol{0}
\end{pmatrix}
+ \boldsymbol C_n^\tau \otimes \begin{pmatrix}
  \boldsymbol M_h\\[3pt] \boldsymbol 0
\end{pmatrix} \boldsymbol V_{n-1}
\,.
\end{equation}
\end{problem}

We note that the orthogonality condition for the pressure has not been implemented yet in Problem~\ref{Prob:LocAlg}. This will be done in the multigrid method by applying the occuring operators to a subspace of $\boldsymbol R^{(k+1)\cdot M^p}$, introduced in~\eqref{Eq:DefRsSs} below. For~\eqref{Eq:InAlg} along with~\eqref{Eq:DefKMC} and~\eqref{Eq:DefMAB} we introduce the abbreviations
\begin{equation}
\label{Eq:DefDtauh}
\boldsymbol D^n_{\tau,h}  \coloneq \begin{pmatrix}
  \boldsymbol K_n^\tau \otimes \boldsymbol M_h + \boldsymbol M_n^\tau \otimes \boldsymbol A_h  &  \boldsymbol M_n^\tau \otimes \boldsymbol B_h^\top \\[3pt]
  \boldsymbol M_n^\tau \otimes \boldsymbol B_h & \boldsymbol 0
\end{pmatrix}\,, \quad
\boldsymbol C^n_{\tau,h}  \coloneq -  \boldsymbol C_n^\tau \otimes \begin{pmatrix}
\boldsymbol M_h\\[3pt] \boldsymbol 0
\end{pmatrix}\,.
\end{equation}
From the local system~\eqref{Eq:InAlg} on $I_n$ we then get the following global problem on $I$.
\begin{problem}[Global algebraic problem]
\label{Prob:GloAlg}
Let $\boldsymbol V_0$ be defined by~\eqref{Eq:DefVnm1}. Find $\boldsymbol{X}=(\boldsymbol{X}_1,\ldots, \boldsymbol{X}_N)\in \R^{N\cdot (k+1) \cdot (M^{\boldsymbol{v}}+M^p)}$ such that
\begin{equation}
\label{Eq:GloSys}
\begin{pmatrix}
\boldsymbol D^1_{\tau,h}  \\[3pt]
\boldsymbol C^2_{\tau,h}  & \boldsymbol D^2_{\tau,h} \\[3pt]
& \ddots & \ddots \\[3pt]
&& \boldsymbol C^N_{\tau,h}  & \boldsymbol D^N_{\tau,h}
\end{pmatrix}
\begin{pmatrix}
  \boldsymbol{X}_1\\[3pt] \vdots \\[3pt] \vdots \\[3pt] \boldsymbol{X}_N
\end{pmatrix}
=
\begin{pmatrix}
\boldsymbol{B}_1 -\boldsymbol C^1_{\tau,h} \boldsymbol{V}_0
\\[3pt] \vdots \\[3pt] \vdots \\[3pt] \boldsymbol{B}_N
\end{pmatrix}\,.
\end{equation}
The subvectors and -matrices in~\eqref{Eq:GloSys} are defined by~\eqref{Eq:DefDtauh} along with~\eqref{Eq:DefKMC} to~\eqref{Eq:DefFn}.
\end{problem}
\begin{remark}[Global linear system of Problem~\ref{Prob:GloAlg}]\label{rem:GLSP}
  \begin{itemize}[leftmargin=*,itemsep=0pt,topsep=0pt,parsep=0pt]
  \item Our $hp$ STMG employs the global system
    representation~\eqref{Eq:GloSys}. The formulation, that orders the unknowns
    by time and local variables as defined in~\eqref{Eq:DefX}, offers the
    advantage of allowing for the restriction of the temporal multigrid to a
    smaller number $\widetilde N< N$ of subintervals by combining $\widetilde N$
    subintervals to a macro time step. This flexibility allows adaptation to the
    available hardware. Conversely, solving the entire system~\eqref{Eq:GloSys}
    for a high number of subintervals demands large computing and memory
    resources, particularly in three space dimensions. For $\widetilde N=1$, a
    time marching scheme is obtained from~\eqref{Eq:GloSys} with algebraic
    system
  \begin{equation}
    \label{Eq:AlgSysTMS}
    \boldsymbol D^n_{\tau,h}  \boldsymbol X_n = \boldsymbol B_n - \boldsymbol C^n_{\tau,h} \boldsymbol X_{n-1}\,,
  \end{equation}
  for $n=2,\ldots,N$ and  right-hand side $\boldsymbol{B}_1 -\boldsymbol
  C^1_{\tau,h} \boldsymbol{V}_0$ for $n=1$. Our implementation supports a
  flexible choice of $\widetilde N$, as shown
  in~\cite{margenbergSpaceTimeMultigridMethod2024a}. As the systems already
  become quite large for small $\widetilde N$, particularly in 3D, we restrict
  ourselves to $\widetilde N=1$ in the numerical experiments in \Cref{sec:experiments}.

\item Alternatively, a global variable and time major order of the unknowns can be applied, such that instead of~\eqref{Eq:DefX} the vector of all unknowns is defined by
\begin{equation}
  \label{Eq:DefXg}
  \boldsymbol  X = (\boldsymbol  X^{\boldsymbol v},\boldsymbol X^p)^\top\coloneq \big(\boldsymbol  V_{1},\ldots, \boldsymbol  V_{N},\boldsymbol P_{1},\ldots, \boldsymbol P_{N}\big)^\top\in \R^{N\cdot (k+1)\cdot ({M^{\boldsymbol v}}  + M^p)}\,,
\end{equation}
with $\boldsymbol V_{n}$ and $\boldsymbol P_{n}$ of~\eqref{Eq:DefSubvec_2}. This
global in time formulation leads to a system matrix with saddle point structure such that block solver techniques, such as Schur complement methods, become feasible. However, the global space-time system  comprises a large
number of unknowns, necessitating substantial computational resources,
particularly in three space dimensions. The approach is not investigated in the present study as well.
\end{itemize}
\end{remark}

\section{\label{sec:mg-framework}Multigrid framework}
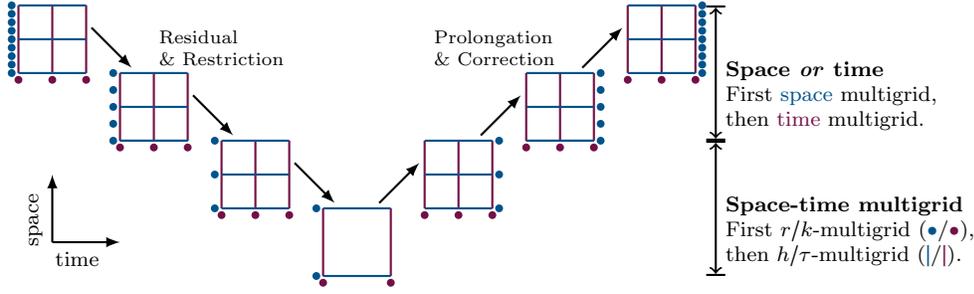
\begin{figure}
  \centering
  \begin{tikzpicture}[xscale=0.9,yscale=0.9,anchor=north west,thick,font=\footnotesize]
    \draw[-latex,shorten >=.1cm,shorten <=.1cm] (.75,1) -- ++(.75,.75);%
    \draw[latex-,shorten >=.1cm,shorten <=.1cm] (.25,1) -- ++(-.75,.75);%
    \draw[shift={(2.5,2)},latex-,shorten >=.1cm,shorten <=.1cm] (.5,.75) --
    ++(-.75,-.75);%
    \draw[shift={(-2.5,2)},-latex,shorten >=.1cm,shorten <=.1cm] (.5,.75) --
    ++(.75,-.75);%
    \draw[shift={(4,3)},latex-,shorten >=.1cm,shorten <=.1cm] (.5,.75) --
    node[left=.5cm,align=left,font=\scriptsize]{Prolongation\\ \& Correction} ++(-.75,-.75);%
    \draw[shift={(-4,3)},-latex,shorten >=.1cm,shorten <=.1cm] (.5,.75) --
    node[right=.5cm,align=left,font=\scriptsize]{Residual\\ \& Restriction} ++(.75,-.75);%
    \foreach \i in {0,1.0}{%
      \draw[myred] (\i, 0) -- ++ (0,1);%
    }%
    \foreach \i in {0,1.0}{%
      \draw[myblue] (0,\i) -- ++ (1,0);%
    }%
    \foreach \i in {0,1}{%
      \filldraw[myred] (\i,-0.1) circle (1.2pt);%
    }%
    \foreach \i in {0,1}{%
      \filldraw[myblue] (-0.1,\i) circle (1.2pt);%
    }%
    \def\points{(-1.5,1), (1.5,1)}; %
    \foreach \p in \points {%
      \foreach \i in {0,0.5,1.0}{%
        \draw[myred,shift={\p}] (\i, 0) -- ++ (0,1);%
      }%
      \foreach \i in {0,0.5,1.0}{%
        \draw[myblue,shift={\p}] (0,\i) -- ++ (1,0);%
      }%
    }%
    \foreach \i in {0,0.5,1}{%
      \filldraw[myred,shift={(-1.5,1)}] (\i,-0.1) circle (1.2pt);%
    }%
    \foreach \i in {0,0.5,...,1}{%
      \filldraw[myblue,shift={(-1.5,1)}] (-0.1,\i) circle (1.2pt);%
    }%
    \foreach \i in {0,0.5,1}{%
      \filldraw[myred,shift={(1.5,1)}] (\i,-0.1) circle (1.2pt);%
    }%
    \foreach \i in {0,0.5,...,1}{%
      \filldraw[myblue,shift={(1.5,1)}] (1.1,\i) circle (1.2pt);%
    }%
    \def\points{(3,2), (-3,2)}%
    \foreach \p in \points {%
      \foreach \i in {0,0.5,1.0}{%
        \draw[myred,shift={\p}] (\i, 0) -- ++ (0,1);%
      }%
      \foreach \i in {0,0, 0.5,1.0}{%
        \draw[myblue,shift={\p}] (0,\i) -- ++ (1,0);%
      }%
    }%
    \foreach \i in {0,0.5,1}{%
      \filldraw[myred,shift={(-3,2)}] (\i,-0.1) circle (1.2pt);%
    }%
    \foreach \i in {0,0.25,...,1}{%
      \filldraw[myblue,shift={(-3,2)}] (-0.1,\i) circle (1.2pt);%
    }%
    \foreach \i in {0,0.5,1}{%
      \filldraw[myred,shift={(3,2)}] (\i,-0.1) circle (1.2pt);%
    }%
    \foreach \i in {0,0.25,...,1}{%
      \filldraw[myblue,shift={(3,2)}] (1.1,\i) circle (1.2pt);%
    }%
    \def\points{(4.5,3), (-4.5,3)}%
    \foreach \p in \points {%
      \foreach \i in {0,0.5,1.0}{%
        \draw[myred,shift={\p}] (\i, 0) -- ++ (0,1);%
      }%
      \foreach \i in {0,0.5,1}{%
        \draw[myblue,shift={\p}] (0,\i) -- ++ (1,0);%
      }%
    }%
    \foreach \i in {0,0.5,1}{%
      \filldraw[myred,shift={(-4.5,3)}] (\i,-0.1) circle (1.2pt);%
    }%
    \foreach \i in {0,0.125,...,1}{%
      \filldraw[myblue,shift={(-4.5,3)}] (-0.1,\i) circle (1.2pt);%
    }%
    \foreach \i in {0,0.5,1}{%
      \filldraw[myred,shift={(4.5,3)}] (\i,-0.1) circle (1.2pt);%
    }%
    \foreach \i in {0,0.125,...,1}{%
      \filldraw[myblue,shift={(4.5,3)}] (1.1,\i) circle (1.2pt);%
    }%
    \draw[-latex] (-4,0.5) -- node[pos=-.1]{time} ++ (1,0);%
    \draw[-latex] (-4,0.5) -- node[rotate=90,above,pos=.33]{space} ++ (0,1);%
    \draw[|<->|] (5.8,0) -- node[align=left,pos=.66]{{\bfseries Space-time multigrid}\\
      First $r/k$-multigrid ({\color{myblue}$\bullet$}/{\color{myred}$\bullet$}),\\
      then $h/\tau$-multigrid ({\color{myblue}$\boldsymbol\vert$}/{\color{myred}$\boldsymbol\vert$}).} ++ (0,2);%
    \draw[|<->|] (5.8,2) -- node[align=left,pos=.66]{{\bfseries Space \emph{or} time}\\
      First {\color{myblue} space} multigrid,\\
      then {\color{myred} time} multigrid.} ++ (0,2);%
  \end{tikzpicture}
  \caption{Sketch of the $hp$ STMG of \Cref{alg:stmg}. The
    corrections are transferred by the prolongation operators and the residual
    is transferred by the restriction operators. On each level the error is
    smoothened by application of the Vanka operator~\eqref{vanka0}.
    The coarsening strategy which is
    used in \Cref{alg:CombineHierarchies}, is first in space and then in
    time in combination with polynomial coarsening before geometric coarsening (cf.~\ref{itm:p1},~\ref{itm:p2}).
  }\label{fig:stmg}
\end{figure}
For solving~\eqref{Eq:GloSys} efficiently, we propose an $hp$ space-time
multigrid method. Here, we use it as a preconditioner for GMRES iterations
rather than as a solver itself. For foundational principles of multigrid methods
we refer
to~\cite{hackbuschMultigridMethodsApplications1985,brambleMultigridMethods1993,vassilevskiMultilevelBlockFactorization2008}.
Before we present the $hp$ STMG in \Cref{alg:stmg}, that is sketched in
\Cref{fig:stmg}, we need to define the grid transfer operators, restriction and
prolongation, and the smoother.

 Let $\{\mathcal{M}_l\}_{l=0}^{L}$ be a quasi-uniform family of nested triangulations of the interval $I$ into semi-closed subintervals $(t_a,\,t_b]$ based on \textit{global regular refinement}, with $\mathcal M_{l} = \{I_i=(t_{i,a},\,t_{i,b}] \mid i=1,\ldots , N^{\text{el}}_{l}\}$, for $l=0,\ldots,L$. The finest partition of $I$ is $\mathcal M_\tau=\mathcal M_{L}$. For the characteristic mesh size $\tau_l$ there holds $\tau_l= \frac{\tau_{l-1}}{2}$, and $\tau_0 = \mathcal O(1)$. This results in a hierarchy of nested temporal finite element spaces  of type~\eqref{Eq:BasYk},
 \begin{equation}
 \label{Eq:DefYkl}
 Y^k_0 \subset  Y^k_1  \subset \cdots \subset Y^k_L \subset L^2(0,\,T) \,.
 \end{equation}
Let further $\{\mathcal{T}_{s}\}_{s=0}^{S}$ be a quasi-uniform family of nested triangulations of the spatial domain $\Omega$ into (open) quadrilaterals or hexahedrals based on \textit{global regular refinement}, with $\mathcal T_{s} = \{K_i\mid i=1,\ldots , N^{\text{el}}_{s}\}$, for $s =0,\ldots,S$. The finest partition is $\mathcal T_h=\mathcal T_{S}$. For the characteristic mesh size $h_s$ there holds $h_s \approx \frac{h_{s-1}}{2}$ and $h_0 = \mathcal O(1)$. This results in a hierarchy of nested spatial finite element spaces of type~\eqref{Eq:BasVhQh},
 \begin{subequations}
 \label{Eq:NFES}
 \begin{alignat}{6}
 \label{Eq:NFESv}
\boldsymbol V^{r+1}_0(\Omega) &  \subset \boldsymbol V^{r+1}_1(\Omega) & \subset \cdots & \subset \boldsymbol V^{r+1}_S(\Omega) && \subset \boldsymbol V\,,\\
 \label{Eq:NFESp}
Q^r_0(\Omega)&  \subset  Q^r_1(\Omega) & \subset \cdots & \subset  Q^r_S(\Omega) && \subset  Q\,.
 \end{alignat}
 \end{subequations}
Similarly to~\eqref{Def:Qhp}, for the nested finite element spaces~\eqref{Eq:NFESp} we define
\begin{equation}
\label{Eq:NFESp+}
Q_s^{r,+} \coloneq Q_s^{r}+\operatorname{span}\{1\}\,.
\end{equation}
We put $M^{\boldsymbol v}_{r+1,s} \coloneq \operatorname{dim}\; \boldsymbol V_s^{r+1}$ and $M^{p}_{r,s} \coloneq \operatorname{dim}\; Q_s^{r,+}$. On the space-time and polynomial order multigrid hierarchy, we define the algebraic tensor product spaces, for $l=0,\ldots,L$ and $s=0,\ldots,S$, as well as $k,r\in \N$ by
\begin{equation}
\label{Eq:STSMG}
\boldsymbol H_{l,s}^{k,r+1}  \coloneq  Y_l^{k}(I) \otimes \boldsymbol V_s^{r+1}(\Omega)\,,\qquad
H_{l,s}^{k,r}  \coloneq Y_l^{k}(I)  \otimes Q_s^{r}(\Omega)\,.
\end{equation}
All spaces, quantities and problems that were introduced before in
\Cref{sec:stfem} and \Cref{sec:AlgSys} are now studied on the full
space-time multigrid hierarchy and for the function spaces~\eqref{Eq:STSMG}. For this, the respective indices are added in the notation. The grid levels in time and space are denoted by the indices ``$l$'' and ``$s$'', respectively, and the polynomial orders by the indices ``$k$'' and ``$r$''.

\subsection{\label{Sec:GTO}Grid transfer operators}
For the multigrid hierarchy of discrete spaces~\eqref{Eq:STSMG}, we need to
define suitable transfer operators.
For our $hp$-STMG method, we need to define them for spatio-temporal geometric refinement
and coarsening, referred to as $h$-multigrid, and for spatio-temporal mesh refinement
and coarsening of the polynomial degree \(k,r\in\mathbb{N}\), referred to as $p$-multigrid.
In both cases, we use the canonical choices of the restriction and prolongation
operators regarding the hierarchy of the discrete spaces~\eqref{Eq:STSMG}. In
Definition~\ref{Def:RP} the notation for the action of the transfer operators
is introduced. For the tighter and technical definition of the restriction and
prolongation operators we refer to Appendix~\ref{App:TOMG}.

\begin{definition}[Restriction and prolongation]
\label{Def:RP}
Let  $k = 2^K$ and $r = 2^R$ for some $K,R \in \mathbb{N}$. Let $l\in \{0,\ldots,L\}$ and $s\in \{0,\ldots,S\}$. Let $D^{k,r}_{l,s}$ denote the total number of degrees of freedom associated with the product space $\boldsymbol H_{l,s}^{k,r+1} \times  Q_s^{r}(\Omega)$; cf.~\eqref{Eq:STSMG}. 
\begin{itemize}[leftmargin=*,itemsep=0pt,topsep=0pt,parsep=0pt]
\item[i)] The prolongation and restriction operators for the space-time mesh are defined as mappings
\begin{equation}
\label{Eq:DefPRSTM}
T_{l-1,l;s-1,s}^{k,r}: \mathbb R^{D^{k,r}_{l-1,s-1}}\to \mathbb R^{D^{k,r}_{l,s}} \quad \text{and} \quad \boldsymbol T^{k,r}_{l,l-1;s,s-1}: \mathbb R^{D^{k,r}_{l,s}}\to \mathbb R^{D^{k,r}_{l-1,s-1}}\,.
\end{equation}

\item[ii)] The prolongation and restriction operators for the polynomial degrees are defined as mappings
\begin{equation}
\label{Eq:DefPRPO}
\boldsymbol T_{l,s}^{\frac{k}{2},k;\frac{r}{2}r,r} : \mathbb R^{D^{\frac{k}{2},\frac{r}{2}}_{l,s}}\to \mathbb R^{D^{k,r}_{l,s}}   \quad \text{and} \qquad
\boldsymbol T_{l,s}^{k,\frac{k}{2};r,\frac{r}{2}}:  \mathbb R^{D^{k,r}_{l,s}} \to   \mathbb R^{D^{\frac{k}{2},\frac{r}{2}}_{l,s}} \,.
\end{equation}

\end{itemize}
\end{definition}
In~\eqref{Eq:DefPRSTM}, the lower touples of indices indicate the direction of changes in the mesh hierarchy. In~\eqref{Eq:DefPRPO}, the upper pairs of indices show the change in the polynomial orders. For the assumption that $k = 2^K$ and $r = 2^R$ we refer to Appendix~\ref{App:TOMG}. 

\subsection{Space-time Vanka smoother}
The smoother for appoximating the solutions of the linear systems on the
  multigrid levels is a key component of multigrid techniques. It aims at smoothing out high frequency errors in the solutions to the linear systems in the $hp$ multigrid hierarchy. In the sequel, these linear systems are assumed to be defined by 
\begin{equation}
  \label{Eq:LSMG}
  \boldsymbol S^{k,r}_{l,s} \boldsymbol X^{k,r}_{l,s}  = \boldsymbol B^{k,r}_{l,s}\,.
\end{equation}
On the finest $hp$ STMG level,~\eqref{Eq:LSMG} represents the linear
system of Problem~\ref{Prob:GloAlg}, or~\ref{Prob:LocAlg} in a time marching
approach. On the coarser levels, the right hand side in~\eqref{Eq:LSMG}
corresponds to the residual, and the solution yields the correction.
The \(V\)-cycle $hp$ STMG, introduced in \Cref{sec:alg}, uses a simple iteration as its smoothing operation,
\[
{\mathtt{smoother}}\Big(\boldsymbol S^{k,r}_{l,s},\,\boldsymbol B^{k,r}_{l,s}\Big) \approx \Big( \boldsymbol S^{k,r}_{l,s}\Big)^{-1}\boldsymbol B^{k,r}_{l,s}\,,
\]
and reduces the high frequency components of the residual
\(\boldsymbol B^{k,r}_{l,s} -\boldsymbol S^{k,r}_{l,s}\boldsymbol
X^{k,r}_{l,s}\). The linear complexity of the geometric multigrid method is
achieved when the reduction rate of the residual is constant across levels. Here, we use a space-time Vanka smoother. For a detailed
description within STFEMs we refer
to~\cite{anselmannGeometricMultigridMethod2023,anselmannEnergyefficientGMRESMultigrid2024}.
For discontinuous pressure, cell-based patches can yield mesh-independent
convergence, as the velocity-pressure coupling is resolved locally~\cite{john_higher-order_2001}. For related multigrid analysis of additive
Schwarz-type smoothers for saddle-point problems, see~\cite{schoberl_schwarz-type_2003}.
For the local (time-stepping) approach~\eqref{Eq:AlgSysTMS}, the cell based Vanka smoother
is built for all $(k+1)\cdot (M^{\boldsymbol v}_{r+1,s}+M^p_{r,s})$ degrees of
freedom of a space-time element which amounts to a block size of
$(k+1)(d{(r+2)}^d+(r+1)^d)$, with space dimension $d$. We use an inner direct
solver. For macro time steps and the global in time approach (cf.\
Remark~\ref{rem:GLSP}) the Vanka smoother is assembled over the subintervals,
which increases the block size. More precisely, on the space-time mesh
\( \mathcal \mathcal T_h \otimes \mathcal M_\tau \) the Vanka smoother is
defined, for \(\boldsymbol S\coloneq \boldsymbol S^{k,r}_{l,s}\) and
\(\boldsymbol b\coloneq \boldsymbol B^{k,r}_{l,s}\), by
\begin{equation}\label{vanka0}
  {\mathtt{smoother}}(\boldsymbol S,\,\boldsymbol b)
  =\bigg(\sum_{T\in \mathcal T_h \otimes \mathcal M_\tau}\boldsymbol R_{T}^\top {[\boldsymbol R_T \boldsymbol S_T\boldsymbol
    R_T^\top]}^{-1}\boldsymbol R_T\bigg) \boldsymbol b\,,
\end{equation}
where \(\boldsymbol R_{T}\) is the restriction to those nodes that belong to the
space-time mesh element \(T\in\mathcal \mathcal T_h \otimes \mathcal M_\tau\)
and $\boldsymbol S_T$ is the corresponding local system matrix on $T$. The
smoother is computationally expensive, for cell based patches its application has
a complexity of \(O\big({(k+1)}^2{(d{(r+2)}^d+(r+1)^d)}^2\big) \). With a relaxation
parameter \(\omega_\ell\in (0,\,1)\), a smoother iteration
(cf.~\Cref{alg:stmg}) is then given by
\begin{equation}\label{vanka}
  \mathtt{Smoother}(\boldsymbol S,\boldsymbol b,\boldsymbol u)=\boldsymbol u + \omega \,{\mathtt{smoother}}(\boldsymbol S,\,\boldsymbol b -\boldsymbol S \boldsymbol u).
\end{equation}
  \begin{remark}[Construction of cell and vertex star patches]
    For the cell based Vanka smoother, we use the standard additive scheme,
    collecting all dofs on each cell and weighting them by their valence. For
    the vertex star patches (VSPs), we follow the construction principle
    in~\cite{anselmannEnergyefficientGMRESMultigrid2024,rafiei2025improvedmonolithicmultigridmethods},
    selecting all velocity dofs and the pressure dofs in the interior of the
    patch. For the $\mathbb{P}_{r}^{\mathrm{disc}}$ pressure used in our
    discretization and in~\cite{anselmannEnergyefficientGMRESMultigrid2024}, all
    pressure dofs are interior and thus included, whereas for continuous
    (e.\,g.\ the
    Taylor-Hood~\cite{rafiei2025improvedmonolithicmultigridmethods}) pressure
    dofs on the boundary of the patch are excluded. This difference arises
    solely from the choice of the pressure space. The velocity dofs are treated
    analogous to the cell based Vanka smoother. In this work, VSPs are applied only in space, with cell-based
    patches in time.
\end{remark}
Our current smoother implementation does not yet exploit the tensor-product
STFEM block structure beyond avoiding assembly of the space-time matrix.
Leveraging it to improve efficiency is left for future work.

\subsection{\label{sec:alg}The {\boldmath$hp$} space-time multigrid method}
\begin{algorithm}[htb]
  \caption{\( hp \) space-time multigrid algorithm for Problem~\ref{Prob:GloAlg}}
  \label{alg:stmg}
  \textbf{Input:} system \( \boldsymbol{S}^{k,r}_{l,s} \boldsymbol{X}^{k,r}_{l,s} = \boldsymbol{B}^{k,r}_{l,s} \); levels \( 0 \le l \le L \), \( 0 \le s \le S \); polynomial orders \( 1\leq k \leq k_{\max} =2^K \), \(1 \leq r \leq r_{\max}=2^R\), $K\leq R$; smoother \( \boldsymbol{W}^{k,r}_{l,s} \); smoothing steps \( \nu_1, \nu_2 \).
  \begin{enumerate}[leftmargin=*]\itemsep1pt \parskip0pt \parsep0pt
    \item[\textbf{1.}] \textbf{Presmoothing.} With an initial guess \( \boldsymbol{X}^{k,r;0}_{l,s} \in \mathbb{R}^{D_{l,s}^{k,r}} \), apply \( \nu_1 \) smoothing steps:
    \begin{equation}\label{Alg:VGMGPre}
      \boldsymbol{X}^{k,r;\nu+1}_{l,s} = \boldsymbol{X}^{k,r;\nu}_{l,s} - (\boldsymbol{W}^{k,r}_{l,s})^{-1} \left( \boldsymbol{S}^{k,r}_{l,s} \boldsymbol{X}^{k,r;\nu}_{l,s} - \boldsymbol{B}^{k,r}_{l,s} \right), \quad \nu = 0, \ldots, \nu_1 - 1.
    \end{equation}
  \item[\textbf{2.}] \textbf{Coarse-grid correction in \( p \).}\hfill$\blacktriangleright$~$p$-Multigrid before $h$-Multigrid.\\
    \textbf{If} \MB{\( r > 1 \)}, compute:
     \[
      \boldsymbol{Y}_{{p}}^{k',\frac{r}{2}} = \textsc{P-Residual-Restrict}( \boldsymbol{X}^{k,r;\nu_1}_{l,s}, \boldsymbol{B}^{k,r}_{l,s}, \boldsymbol{S}^{k,r}_{l,s} )
      \]
      with \( k' = k \) or \( k/2 \), depending on the
     hierarchy.\hfill$\blacktriangleright$~Depends on $ R> K$ or $R=K$;
     
    \textbf{Else}, set: \( \boldsymbol{Y}_{{p}}^{1,1} = \boldsymbol{0} \)\hfill$\blacktriangleright$~Lowest $p$-multigrid level reached; 
  \item[\textbf{3.}] \textbf{Coarse-grid correction in \( h \).} \\
    \textbf{If} \( s > 1 \), compute:
    \[
    \boldsymbol{Y}^{{h}}_{l',s-1} = \textsc{H-Residual-Restrict}( \boldsymbol{X}^{1,1;\nu_1}_{l,s}, \boldsymbol{B}^{1,1}_{l,s}, \boldsymbol{S}^{1,1}_{l,s} )
    \]
    with \( l' = l \) or \( l-1 \), depending on the
    hierarchy\hfill$\blacktriangleright$~Depends on $S>L$ or $S=L$;
    
    \textbf{Else If}  \( s = 1 \), solve: \hfill$\blacktriangleright$~Direct solve on coarse level;
    \[
      \boldsymbol{Y}^{{h}}_{l-1,s-1} \in \mathbb{R}^{D_{l-1,s-1}^{1,1}} \text{ s.\,t.\ } \boldsymbol{S}^{1,1}_{l-1,s-1} \boldsymbol{Y}^{{h}}_{l-1,s-1} = \boldsymbol{B}^{1,1}_{l-1,s-1}
    \] 
  \item[\textbf{4.}] \textbf{Prolongation and correction.} \\
    \textbf{If} \( r = 1 \), compute:
    \[
      \boldsymbol{X}^{1,1;\nu_1+1}_{l,s} = \textsc{H-Prolongate-Correct}( \boldsymbol{X}^{1,1;\nu_1}_{l,s}, \boldsymbol{Y}^{{h}}_{l',s-1} )
    \]
    \textbf{Else}, compute:
    \[
      \boldsymbol{X}^{k,r;\nu_1+1}_{l,s} = \textsc{P-Prolongate-Correct}( \boldsymbol{X}^{k,r;\nu_1}_{l,s}, \boldsymbol{Y}_{{p}}^{k',\frac{r}{2}} )
    \]
  \item[\textbf{5.}] \textbf{Postsmoothing.} Apply \( \nu_2 \) smoothing steps to get:
    \begin{equation}\label{Alg:VGMGPost}
      \boldsymbol{X}^{k,r;\nu+1}_{l,s} = \boldsymbol{X}^{k,r;\nu}_{l,s} - (\boldsymbol{W}^{k,r}_{l,s})^{-1} \left( \boldsymbol{S}^{k,r}_{l,s} \boldsymbol{X}^{k,r;\nu}_{l,s} - \boldsymbol{B}^{k,r}_{l,s} \right), \quad \nu = \nu_1 + 1, \ldots, \nu_1 + \nu_2
    \end{equation}
    \item[\textbf{6.}] \textbf{Return} \( \boldsymbol{X}^{k,r;\nu_1+\nu_2+1}_{l,s} \)
  \end{enumerate}
\end{algorithm}
We introduce our $V$-cycle $hp$ space-time multigrid (STMG) method for tensor-product STFEMs, assuming the linear system of Problem~\ref{Prob:GloAlg} is represented on multigrid levels $l=0,\ldots,L$ (time) and $s=0,\ldots,S$ (space) with polynomial degrees $k,r\in\mathbb{N}$, cf.~\eqref{Eq:LSMG}. 
Although $k=0$ or $r=0$ is theoretically possible, these cases are excluded due
to poor performance in practice. 
For clarity, and without loss of generality, we assume $S \geq L$ and $r \geq
k$, which is typical in FEM flow simulations dominated by spatial dynamics. 
The $hp$~STMG algorithm for~\eqref{Eq:LSMG} is outlined in \Cref{alg:stmg}, with its subroutines in \Cref{alg:pmg}--\ref{alg:ppro} and illustrated in \Cref{fig:stmg}. 
  We assume $k=2^K$ and $r=2^R$ with $K,R\in\mathbb{N}$, $K\leq
  R$; see Appendix~\ref{App:TOMG}. We note that the $hp$~STMG is used as a preconditioner for GMRES.
\begin{remark}[Choice of the iterative solver for~\eqref{Eq:LSMG}]\label{rem:iterative-choice}
  Unlike in the stationary setting, the space-time systems arising in our
  formulation are non-symmetric, even at the level of the local
  problem~\ref{Eq:InAlg}, due to the contribution of the time derivative and the
  resulting non-symmetric matrix \(\boldsymbol{K}_n^\tau\). Consequently,
  some solvers for the stationary Stokes, such as MINRES, are not applicable.
  While block preconditioners may help recover symmetry at the block level, we
  consider such approaches beyond the scope of this work. We therefore use
  GMRES, which, though more memory-intensive, proves efficient and exhibits
  robust convergence behavior in our experiments.
\end{remark}

\subsection{Multigrid Sequence Generation}
\begin{algorithm}[htbp]
  \caption{\textsc{ConstructHierarchy}$(h_L,\,h_0,\,r_L,\,r_0)$}\label{alg:ConstructHierarchy}
  \begin{algorithmic}[1]
    \REQUIRE{Fine and coarse mesh size $h_{L}$, $h_{0}$, polynomial degrees $r_{L}$, $r_{0}$.}
    \STATE{$\mathcal{H} \gets (),\;\;h \gets h_{L_{g}},\;\;r \gets r_{L}$}
    \WHILE{$r\geq r_0$}
    \STATE{$\mathcal{H} \gets (h,\,r) \times \mathcal{H}$} \COMMENT{Polynomial coarsening}
    \STATE{$r \gets \lfloor r/2 \rfloor$}  \COMMENT{halve polynomial degree}
    \ENDWHILE{}
    \WHILE{$h\leq h_0$}
    \STATE{ $\mathcal{H} \gets (h,\,r) \times \mathcal{H}$}\COMMENT{Geometric coarsening}
    \STATE{ $h \gets 2 h$}     \COMMENT{Coarsen triangulation (``double mesh
      size'' for simplicity)}
    \ENDWHILE{}
    \RETURN{$\mathcal{H}$}
  \end{algorithmic}
\end{algorithm}
\begin{algorithm}[htbp]
  \caption{\textsc{CombineHierarchies}$(\mathcal{H}_h,\,\mathcal{H}_\tau)$}\label{alg:CombineHierarchies}
\begin{algorithmic}[1]
  \REQUIRE{ Spatial hierarchy obtained by \Cref{alg:ConstructHierarchy}: $\mathcal{H}_h=\left((h_\ell,\,p_\ell)\right)_{\ell=1}^{\mathcal{S}}$}
  \REQUIRE{} Temporal hierarchy obtained by \Cref{alg:ConstructHierarchy}: $\mathcal{H}_\tau=\left((\tau_\ell,\,k_\ell)\right)_{\ell=1}^{\mathcal{L}}$
\STATE{ $\mathcal{H}_{st} \gets (),\quad(\tau_{\text{pad}},\,k_{\text{pad}}) \gets \mathcal{H}_\tau[\mathcal{L}]$}
\FOR{$\ell = \mathcal{L} + 1$ to $\mathcal{L}_{\max}$}
\STATE{ $\mathcal{H}_\tau \gets \mathcal{H}_\tau \times (\tau_{\text{pad}},\,k_{\text{pad}})$}\label{alg:combinehierarchies:pad}\COMMENT{Pad the hierarchy to length $\mathcal{L}_{\max}$}
\ENDFOR{}
\FOR{$\ell = 1$ to $\mathcal{L}_{\max}$}
\STATE{Let $(\tau_\ell,\,k_\ell) \gets \mathcal{H}_\tau[\ell],\quad (h_\ell,\,r_\ell)
  \gets \mathcal{H}_h[\ell]$}\COMMENT{Combine level by level}
\STATE{ $\boldsymbol{H}_{l_{\ell},s_{\ell}}^{\,k_\ell,\,r_\ell} \coloneq
  Y_{l_\ell}^{\,k_\ell}(I)\otimes\boldsymbol{V}_{s_{\ell}}^{\,r_\ell+1}(\Omega),\quad
  H_{l_{\ell},s_{\ell}}^{\,k_\ell,\,r_\ell} \coloneq
Y_{l_\ell}^{\,k_\ell}(I)\otimes Q_{s_{\ell}}^{\,r_\ell}(\Omega)$}
\STATE{$\mathcal{H}_{st} \gets \mathcal{H}_{st} \times
  (\boldsymbol{H}_{l_{\ell},s_{\ell}}^{\,k_\ell,\,r_\ell},\, H_{l_{\ell},s_{\ell}}^{\,k_\ell,\,r_\ell})$}\label{alg:combinehierarchies:stfespaces}
  \ENDFOR{}
 \RETURN{ $\mathcal{H}_{st}$}
\end{algorithmic}
\end{algorithm}
Finally, we comment on some implementational aspects for Algorithm~\ref{alg:stmg}. We construct the space-time multigrid hierarchy by coarsening sequences of
spatial and temporal finite element spaces, which are combined by tensor
products, according to two guiding principles:
\begin{itemize}[itemsep=0pt,topsep=0pt,parsep=0pt]
\item[\customlabel{itm:p1}{(P1)}]
  \textbf{Spatial Coarsening over Temporal Coarsening}:
  Perform geometric coarsening first in the spatial dimension, then in the temporal dimension.
\item[\customlabel{itm:p2}{(P2)}]
  \textbf{Polynomial over Geometric Coarsening}:
  Apply coarsening in polynomial degrees \((r\)/\(k\)) before geometric coarsening \((h\)/\(\tau\)).
\end{itemize}
Following~\ref{itm:p2}, the sequences of nested temporal and spatial finite
element spaces are generated through \Cref{alg:ConstructHierarchy}. We generate
a temporal hierarchy of finite element spaces \(Y_l^k(I)\)
(cf.~\eqref{Eq:DefYkl}) with geometric level $l$ and polynomial degree $k$. Let
\(L\) be the number of geometric levels in time, and
\(\lfloor \log_2(k)\rfloor\) the number of polynomial levels. Thus, the number
of temporal levels is \(\mathcal{L} \coloneq L + \lfloor \log_2(k)\rfloor\).
The spatial hierarchy of finite element spaces
\(\boldsymbol{V}_s^{r+1}(\Omega),\,Q_s^r(\Omega)\) on geometric level $s$ with
polynomial degree $r$ (cf.~\eqref{Eq:NFES}) is generated analogously. Let \(S\)
be the number of geometric levels in space, and \(\lfloor \log_2(r)\rfloor\) the
number of polynomial levels. Thus, the number of spatial levels is
\(\mathcal{S} \coloneq S + \lfloor \log_2(r)\rfloor\). These two hierarchies
are merged into a single STMG sequence of size
\(\mathcal{L}_{\max} = \max(\mathcal{L}, \mathcal{S})\). This is represented in \Cref{alg:CombineHierarchies}. If
\(\mathcal{L} < \mathcal{L}_{\max}\), we pad the final temporal spaces to match
the finer levels in time, setting
\(Y_m^k(I) = Y_{\mathcal{L}}^k(I), m=L+1,\dots,S\)
(cf. \Cref{alg:CombineHierarchies} line~\ref{alg:combinehierarchies:pad}). If
\(\mathcal{S} < \mathcal{L}_{\max}\), an analogous padding applies to
\(\boldsymbol{V}_m^{r+1}(\Omega)\) and $Q_m^r(\Omega)$ if
\(\mathcal{L}>\mathcal{S}\). For a level
\(\ell\in\{1,\dots,\mathcal{L}_{\max}\}\), we denote $l_\ell$ and $s_\ell$ as
the temporal and spatial geometric level, $k_\ell$ and $r_{\ell}$ as the
temporal and spatial polynomial degree. Following~\ref{itm:p1}, the
space-time finite element spaces
\(\boldsymbol{H}_{l_\ell,s_\ell}^{\,k_\ell,\,r_\ell+1}\) and
\(H_{l_\ell,s_\ell}^{\,k_\ell,\,r_\ell}\), where
\(\ell\in\{1,\dots,\mathcal{L}_{\max}\}\), are generated in
\Cref{alg:CombineHierarchies} line~\ref{alg:combinehierarchies:stfespaces}.
According to the construction in \Cref{alg:CombineHierarchies}
and~\ref{alg:ConstructHierarchy}, the spaces are ordered such that the coarsest
space-time function space, characterized by large $h$, $\tau$ and small $p$,
$k$, is located at $\ell=0$. Two additional principles for the generation of the
multgrid hierarchy naturally arise from \Cref{alg:CombineHierarchies} and~\ref{alg:ConstructHierarchy}.
\begin{itemize}[itemsep=0pt,topsep=0pt,parsep=0pt]
\item[\customlabel{itm:p3}{(P3)}] \textbf{True Space-Time Multigrid at Coarsest
    Level}: Levels with space-time coarsening are put at the lower levels of the hierarchy.
\item[\customlabel{itm:p4}{(P4)}] \textbf{Padding for Pure Space or Time Level}:
  Padding with identical function spaces is done at the finer
  levels.
\end{itemize}
For an example of a generated Multigrid sequence, we refer to Appendix~\ref{sec:example-mgseq}.
Although placing full space-time levels at the top of the multigrid hierarchy
could more rapidly reduce the total number of space-time degrees of freedom,
they are placed according to~\ref{itm:p3}. This choice is guided by a CFL-type
condition that is derived in~\cite{chaudet-dumasOptimizedSpaceTimeMultigrid2023}
for one-dimensional space-time multigrid approaches to the heat equation, which
ensures convergence under space-time coarsening. Numerical experiments indicate
that such a condition also arises for the Stokes system: adopting ``early'' space-time coarsening at the top of the hierarchy can degrade the performance of multigrid methods for STFEMs.
\subsection{Matrix-free operator evaluation}\label{sec:mf}
In this work, linear operators are evaluated without the explicit formation and
storage of system matrices. For this, we rely on the matrix-free multigrid
framework in the \texttt{deal.II}
library~\cite{africa_dealii_2024,kronbichlerGenericInterfaceParallel2012,munchEfficientDistributedMatrixfree2023,fehnHybridMultigridMethods2020}.
A matrix-vector product \(\boldsymbol{Y} =\boldsymbol S \boldsymbol X\)
(cf.~\eqref{Eq:LSMG}) is computed via global accumulation of local element-wise
operations,
\[
  \boldsymbol{S} \boldsymbol X
  =
  \sum_{c=1}^{n_c}
  \boldsymbol{R}_{c,\mathrm{loc\text{-}glob}}^\top\boldsymbol{S}_c\boldsymbol{R}_{c,\mathrm{loc\text{-}glob}}\boldsymbol X ,
  \quad
  \boldsymbol{S}_c
  =
  \boldsymbol{B}_c^\top\boldsymbol{D}_c\boldsymbol{B}_c,
\]
where \(\boldsymbol{R}_{c,\mathrm{loc\text{-}glob}}\) maps local degrees of
freedom to global indices, \(\boldsymbol{B}_c\) contains shape function
gradients, and \(\boldsymbol{D}_c\) encodes quadrature weights and material
coefficients. Sum-factorization then reduces the multi-dimensional operations to
a product of one-dimensional operations. Vectorization
further accelerates the evaluations.

These techniques are used to assemble matrix-vector products with the spatial
operators in~\eqref{Eq:DefMAB}. The temporal matrices in~\eqref{Eq:DefKMC} are
precomputed as described in~\cite{margenbergSpaceTimeMultigridMethod2024a}.
Products of the form
\((\boldsymbol{M}^{\tau} \otimes \boldsymbol{A}_h)\boldsymbol{V}^a\), for $a=1,\ldots,k+1$ (cf.\ \Cref{Sec:Prem} and~\ref{Sec:Afdp}), are evaluated
by computing \(\boldsymbol{A}_h\boldsymbol{V}^a\) once for each temporal degree
of freedom, followed by a matrix-vector multiplication with the temporal
matrices in a blockwise sense
(cf.~\cite{margenbergSpaceTimeMultigridMethod2024a}). This approach extends
naturally to all Kronecker products in~\eqref{Eq:InAlg}.

\section{\label{sec:experiments}Numerical experiments}
We validate the accuracy and convergence properties of the
proposed $hp$ STMG solver for the Stokes system using a sequence of polynomial
degrees and mesh refinements. We use inf-sup stable
$\mathbb{Q}_{r+1}/\mathbb{P}_{r}^{\text{disc}}$ elements in space and a DG$(k)$
discretization in time. To solve the linear systems of equations, we use a GMRES
method (cf.~Remark~\ref{rem:iterative-choice}) with a single $V$-cycle $hp$ STMG preconditing step per
iteration. To ensure efficiency and scalability, the number of iterations until
convergence is reached must remain bounded as the mesh size $h$ is reduced or
the polynomial degree $r$ is increased. Thus, we characterize the solver's
performance in terms of:
\begin{description}[style=unboxed,leftmargin=0cm]\itemsep1pt
\item[\textbf{$\boldsymbol h$-robustness}:] Iteration counts remain bounded
  independently of the mesh size \(h\).
\item[\textbf{$\boldsymbol p$-robustness}:] Iteration counts remain bounded
  independently of the degree \(p\).
\end{description}
Robustness is important to ensure that the computational cost increases linearly
with the problem size, preserving the computational complexity.
Further, robustness is essential to avoid prohibitive memory consumption,
especially since GMRES is among the most memory-intensive Krylov solvers. While
the method converges fast in our experiments, its performance depends on the
Arnoldi basis size, which may be constrained by available memory. In our
experience, restarting GMRES once this limit is reached tends to degrade
convergence for the instationary Stokes problem. Consequently, for even larger
problems, the inability to enlarge the basis may limit solver robustness. For
the problem sizes considered here, however, GMRES performs very well and poses
no practical limitations.

The tests were performed on an HPC cluster (HSUper at Helmut Schmidt University) with 571 nodes, each with 2 Intel Xeon Platinum 8360Y CPUs and \SI[scientific-notation=false,round-precision=0]{256}{\giga\byte} RAM.\@ The processors
have 36 cores each and the number of MPI processes always match the cores.
As mentioned in \Cref{rem:GLSP}, we restrict ourselves to $\widetilde N=1$, for
computational studies with $\widetilde N>1$ we refere to~\cite{margenbergSpaceTimeMultigridMethod2024a}.
\subsection{Convergence test}
\label{sec:conv-stokes}
\begin{table}[htb]
  \caption{Errors for
    $\mathbb{Q}_{r+1}^2/\mathbb{P}_{r}^{\mathrm{disc}}/\mathrm{DG}(r)$
    discretizations of the Stokes system
    for~\eqref{eq:conv-test-v}.}\label{tab:conv-stokes}
  \vspace*{-2ex}
  \begin{subcaptionblock}{\textwidth}
    \caption{Calculated velocity and pressure errors in the space-time $L^2$-norm with eoc.}
    \setlength{\tabcolsep}{7pt} \centering \footnotesize
    \begin{tabular}{l|llll|llll}
      \toprule
      &\multicolumn{4}{c|}{$r=4$}&\multicolumn{4}{c}{$r=5$}\\
      $h$ & $e^{\boldsymbol v}_{L^2/L^2}$     &eoc&$e^{p}_{L^2/L^2}$&eoc&$e^{\boldsymbol v}_{L^2/L^2}$ &eoc&$e^{p}_{L^2/L^2}$&eoc\\
      \midrule
      ${2}^{-1}$ & \num{1.00279e-04} &    -&\num{1.14907e-03}&    -&\num{2.71058e-05} &   -  &\num{1.01607e-03} &    -\\
      ${2}^{-2} $& \num{2.32706e-06} & 5.43&\num{1.11534e-04}& 3.36&\num{2.28082e-07} & 6.89 &\num{1.40569e-05} & 6.18\\
      ${2}^{-3} $& \num{3.98114e-08} & 5.87&\num{3.58645e-06}& 4.96&\num{1.87701e-09} & 6.92 &\num{2.29909e-07} & 5.93\\
      ${2}^{-4}$& \num{6.39216e-10} & 5.96&\num{1.12642e-07}& 4.99&\num{1.58961e-11} & 6.88 &\num{3.88539e-09} & 5.89\\
      ${2}^{-5}$& \num{1.10823e-11} & 5.85&\num{3.75088e-09}& 4.91&\num{3.52916e-12} & 2.17 &\num{1.17454e-09} & 1.73\\
      \bottomrule
    \end{tabular}
  \end{subcaptionblock}

  \vspace*{-1ex}
  \begin{subcaptionblock}{\textwidth}
    \caption{Calculated velocity errors in the space-time $H^{1}/L^2$-norm and divergence with eoc.}
    \setlength{\tabcolsep}{7pt} \centering\footnotesize
    \begin{tabular}{l|llll|llll}
      \toprule
      &\multicolumn{4}{c|}{$r=4$}&\multicolumn{4}{c}{$r=5$}\\
      $h$ & $e^{\boldsymbol v}_{H^1/L^2}$     &eoc&$e^{\boldsymbol{\nabla \cdot v}}_{L^2/L^2}$&eoc&$e^{\boldsymbol v}_{H^1/L^2}$ &eoc&$e^{\boldsymbol{\nabla \cdot v}}_{L^2/L^2}$&eoc\\
      \midrule
      ${2}^{-1}$ &\num{4.083771e-03}&  - &\num{6.1022e-04} &   - &\num{8.810760e-04}&   - &\num{8.2991e-04}&    - \\ 
      ${2}^{-2}$ &\num{1.446197e-04}&4.82&\num{1.0657e-04} &2.52 &\num{1.472277e-05}&5.90 &\num{1.2840e-05}& 6.01 \\ 
      ${2}^{-3}$ &\num{4.783616e-06}&4.91&\num{3.6748e-06} &4.86 &\num{2.396664e-07}&5.94 &\num{2.1445e-07}& 5.90 \\ 
      ${2}^{-4}$ &\num{1.524543e-07}&4.97&\num{1.1875e-07} &4.95 &\num{3.760674e-09}&5.99 &\num{3.3924e-09}& 5.98 \\ 
      ${2}^{-5}$ &\num{4.796613e-09}&4.99&\num{3.7541e-09} &4.98 &\num{1.484210e-10}&4.66 &\num{1.1179e-10}& 4.92 \\
      \bottomrule
    \end{tabular}
  \end{subcaptionblock}
\end{table}
\begin{table}[htb]
  \centering
  \caption{Number of GMRES iterations until convergence for
    polynomial degrees $r$ and number of refinements $c$ with  \(\mathbb{Q}_{r+1}^2/\mathbb{P}_{r}^{\text{disc}}/\text{DG}(r)\)
    discretization of the Stokes system and for $hp$ STMG (left) and
    geometric multigrid in space method (right).}\label{tab:iter-stokes}
  \vspace*{-2ex}
  \begin{subcaptionblock}{\textwidth}
    \setlength{\tabcolsep}{5.5pt}
    \caption{Results for the Vanka smoother based on cells.}
  \begin{minipage}{0.47\textwidth}
      \centering\scriptsize
      \begin{tabular}{ccccccc}
        \toprule
        \diagbox[innerleftsep=0pt,innerrightsep=0pt, height=1.5em, width=1.5em]{$r$}{$c$}  &  1  &  2  &  3  &  4  &  5  &  6  \\
        \midrule
                               2& {14.0}& 15.0& 15.0& 14.0& 13.0& 10.6\\
                               3& 19.8& {15.9}& 16.0& 15.0& 13.7& 11.0\\
                               4& 27.8& 23.0& {22.9}& 21.9& 19.0& 15.5\\
                               5& 31.0& 26.4& 26.6& {22.8}& 18.7& 14.9\\
                               6& 45.0& 36.1& 36.7& 29.0& {23.1}& 17.2\\
                               7& 50.8& 43.8& 42.8& 32.8& 25.6& {19.6}\\
        \bottomrule
      \end{tabular}
  \end{minipage}
  \hspace{.3cm}
  \begin{minipage}{0.47\textwidth}
    \centering\scriptsize
      \begin{tabular}{ccccccc}
        \toprule
        \diagbox[innerleftsep=0pt,innerrightsep=0pt, height=1.5em, width=1.5em]{$r$}{$c$} &  1  &  2  &  3  &  4  &  5  &  6  \\
        \midrule
                         2& {14.0}& 15.0& 15.0& 14.0& 13.0& 10.6\\
                         3& 19.0& {17.9}& 18.9& 18.3& 16.4& 14.0\\
                         4& 24.0& 26.8& {24.7}& 24.6& 21.4& 18.4\\
                         5& 26.0& 26.4& 28.8& {27.7}& 24.7& 21.9\\
                         6& 35.0& 33.9& 34.6& 30.9& {29.6}& 26.9\\
                         7& 40.0& 38.8& 39.6& 36.7& 34.5& {31.9}\\
        \bottomrule
      \end{tabular}
    \end{minipage}
  \end{subcaptionblock}

  \begin{subcaptionblock}{\textwidth}
    \setlength{\tabcolsep}{2pt}
    \caption{Total number of elements of all local cell patch matrices in the Vanka smoother $N_{\text{sm}}^{\text{el}}$.}
  \begin{minipage}{0.49\textwidth}
      \centering\scriptsize
      \begin{tabular}{@{}c*{6}{c}@{}}
        \toprule
        \diagbox[innerleftsep=1pt,innerrightsep=-1pt, height=1.5em, width=1.25em]{$r$}{$c$} &  1  &  2  &  3  &  4  &  5  &  6  \\
        \midrule
        2 & \num[round-precision=1]{8820} & \num[round-precision=1]{37044} & \num[round-precision=1]{149940} & \num[round-precision=1]{601524} & \num[round-precision=1]{2407860} & \num[round-precision=1]{9633204} \\
        3 & \num[round-precision=1]{60804} & \num[round-precision=1]{280260} & \num[round-precision=1]{1158084} & \num[round-precision=1]{4669380} & \num[round-precision=1]{18714564} & \num[round-precision=1]{74895300} \\
        4 & \num[round-precision=1]{239220} & \num[round-precision=1]{1043316} & \num[round-precision=1]{4259700} & \num[round-precision=1]{17125236} & \num[round-precision=1]{68587380} & \num[round-precision=1]{274435956} \\
        5 & \num[round-precision=1]{817704} & \num[round-precision=1]{3677524} & \num[round-precision=1]{15198724} & \num[round-precision=1]{61283524} & \num[round-precision=1]{245622724} & \num[round-precision=1]{982979524} \\
        6 & \num[round-precision=1]{2099988} & \num[round-precision=1]{9060804} & \num[round-precision=1]{37042308} & \num[round-precision=1]{148968324} & \num[round-precision=1]{596672388} & \num[round-precision=1]{2387488644} \\
        7 & \num[round-precision=1]{5009076} & \num[round-precision=1]{22023540} & \num[round-precision=1]{89111988} & \num[round-precision=1]{357465780} & \num[round-precision=1]{1430880948} & \num[round-precision=1]{5724541620} \\
        \bottomrule
      \end{tabular}
  \end{minipage}
  \hfill
  \begin{minipage}{0.49\textwidth}
    \centering\scriptsize
    \begin{tabular}{@{}c*{6}{c}@{}}
        \toprule
        \diagbox[innerleftsep=1pt,innerrightsep=-1pt, height=1.5em, width=1.25em]{$r$}{$c$} &  1  &  2  &  3  &  4  &  5  &  6  \\
        \midrule
        2 & \num[round-precision=1]{8820} & \num[round-precision=1]{37044} & \num[round-precision=1]{149940} & \num[round-precision=1]{601524} & \num[round-precision=1]{2407860} & \num[round-precision=1]{9633204} \\
        3 & \num[round-precision=1]{64980} & \num[round-precision=1]{272916} & \num[round-precision=1]{1104660} & \num[round-precision=1]{4431636} & \num[round-precision=1]{17739540} & \num[round-precision=1]{70971156} \\
        4 & \num[round-precision=1]{288000} & \num[round-precision=1]{1209600} & \num[round-precision=1]{4896000} & \num[round-precision=1]{19641600} & \num[round-precision=1]{78624000} & \num[round-precision=1]{314553600} \\
        5 & \num[round-precision=1]{946125} & \num[round-precision=1]{3973725} & \num[round-precision=1]{16084125} & \num[round-precision=1]{64525725} & \num[round-precision=1]{258292125} & \num[round-precision=1]{1033357725} \\
        6 & \num[round-precision=1]{2548980} & \num[round-precision=1]{10705716} & \num[round-precision=1]{43332660} & \num[round-precision=1]{173840436} & \num[round-precision=1]{695871540} & \num[round-precision=1]{2783995956} \\
        7 & \num[round-precision=1]{5962320} & \num[round-precision=1]{25041744} & \num[round-precision=1]{101359440} & \num[round-precision=1]{406630224} & \num[round-precision=1]{1627713360} & \num[round-precision=1]{6512045904} \\
        \bottomrule
      \end{tabular}
    \end{minipage}
  \end{subcaptionblock}

  \begin{subcaptionblock}{\textwidth}
    \setlength{\tabcolsep}{5.5pt}
    \caption{Results for the Vanka smoother based on vertex star patches.}
    \begin{minipage}{0.47\textwidth}
      \centering\scriptsize
      \begin{tabular}{ccccccc}
        \toprule
        \diagbox[innerleftsep=0pt,innerrightsep=0pt, height=1.5em, width=1.5em]{$r$}{$c$}  &  1  &  2  &  3  &  4  &  5  &  6  \\
        \midrule
                               2& {12.5}& 15.0& 15.0& 14.0& 13.0& 11.9\\
                               3& 18.0& {15.6}& 15.3& 15.0& 14.9& 12.0\\
                               4& 36.8& 24.5& {22.0}& 22.2& 20.8& 15.8\\
                               5& 32.5& 24.9& 26.0& {24.0}& 20.0& 16.7\\
                               6& 43.8& 33.5& 35.2& 31.3& {24.7}& 18.9\\
                               7& 45.2& 40.5& 41.2& 33.7& 26.9& {20.0}\\
        \bottomrule
      \end{tabular}
  \end{minipage}
  \hspace{.3cm}
  \begin{minipage}{0.47\textwidth}
    \centering\scriptsize
      \begin{tabular}{ccccccc}
        \toprule
        \diagbox[innerleftsep=0pt,innerrightsep=0pt, height=1.5em, width=1.5em]{$r$}{$c$} &  1  &  2  &  3  &  4  &  5  &  6  \\
        \midrule
                         2& {12.5}& 15.0& 15.0& 14.0& 13.0& 11.9\\
                         3& 17.0& {17.0}& 18.0& 18.0& 16.3& 14.8\\
                         4& 24.0& 37.9& {26.4}& 34.0& 20.6& 26.2\\
                         5& 25.0& 25.0& 26.8& {26.0}& 23.8& 20.8\\
                         6& 34.5& 29.4& 30.3& 29.9& 27.5& {23.6}\\
                         7& 40.5& 36.2& 33.1& 30.6& 27.9& {25.4}\\
        \bottomrule
      \end{tabular}
    \end{minipage}
\end{subcaptionblock}
\begin{subcaptionblock}{\textwidth}
    \setlength{\tabcolsep}{2pt}
    \caption{Total number of elements of all local vertex star patch matrices in the smoother $N_{\text{sm}}^{\text{el}}$.
}
  \begin{minipage}{0.49\textwidth}
    \centering\scriptsize
    \begin{tabular}{@{}c*{6}{c}@{}}
        \toprule
        \diagbox[innerleftsep=1pt,innerrightsep=-1pt, height=1.5em, width=1.25em]{$r$}{$c$} &  1  &  2  &  3  &  4  &  5  &  6  \\
        \midrule
        2 & \num[round-precision=1]{26088} & \num[round-precision=1]{144508} & \num[round-precision=1]{674832} & \num[round-precision=1]{2913700} & \num[round-precision=1]{12108600} & \num[round-precision=1]{49371340} \\
        3 & \num[round-precision=1]{262131} & \num[round-precision=1]{1525200} & \num[round-precision=1]{7047933} & \num[round-precision=1]{30089418} & \num[round-precision=1]{124166103} & \num[round-precision=1]{504303972} \\
        4 & \num[round-precision=1]{1057752} & \num[round-precision=1]{5958040} & \num[round-precision=1]{27384312} & \num[round-precision=1]{116764680} & \num[round-precision=1]{481651224} & \num[round-precision=1]{1955984424} \\
        5 & \num[round-precision=1]{4381406} & \num[round-precision=1]{25621858} & \num[round-precision=1]{118663058} & \num[round-precision=1]{506601133} & \num[round-precision=1]{2089920008} & \num[round-precision=1]{8486454883} \\
        6 & \num[round-precision=1]{12382179} & \num[round-precision=1]{70086273} & \num[round-precision=1]{323040957} & \num[round-precision=1]{1378582665} & \num[round-precision=1]{5688291669} & \num[round-precision=1]{23102415393} \\
        7 & \num[round-precision=1]{30736660} & \num[round-precision=1]{176035016} & \num[round-precision=1]{812986800} & \num[round-precision=1]{3439490827} & \num[round-precision=1]{14187763866} & \num[round-precision=1]{57617435977} \\
        \bottomrule
      \end{tabular}
  \end{minipage}
  \hfill
  \begin{minipage}{0.49\textwidth}
    \centering\scriptsize
    \begin{tabular}{@{}c*{6}{c}@{}}
        \toprule
        \diagbox[innerleftsep=1pt,innerrightsep=-1pt, height=1.5em, width=1.25em]{$r$}{$c$} &  1  &  2  &  3  &  4  &  5  &  6  \\
        \midrule
        2 & \num[round-precision=1]{26088} & \num[round-precision=1]{144508} & \num[round-precision=1]{674832} & \num[round-precision=1]{2913700} & \num[round-precision=1]{12108600} & \num[round-precision=1]{49371340} \\
        3 & \num[round-precision=1]{265806} & \num[round-precision=1]{1498473} & \num[round-precision=1]{7060644} & \num[round-precision=1]{30626847} & \num[round-precision=1]{127581210} & \num[round-precision=1]{520832277} \\
        4 & \num[round-precision=1]{1154176} & \num[round-precision=1]{6612448} & \num[round-precision=1]{31375696} & \num[round-precision=1]{136563840} & \num[round-precision=1]{569858928} & \num[round-precision=1]{2328436256} \\
        5 & \num[round-precision=1]{4605800} & \num[round-precision=1]{26511050} & \num[round-precision=1]{126116575} & \num[round-precision=1]{549656050} & \num[round-precision=1]{2295190275} & \num[round-precision=1]{9381320850} \\
        6 & \num[round-precision=1]{13538628} & \num[round-precision=1]{78401700} & \num[round-precision=1]{374128524} & \num[round-precision=1]{1633231116} & \num[round-precision=1]{6825735972} & \num[round-precision=1]{27911813076} \\
        7 & \num[round-precision=1]{33159133} & \num[round-precision=1]{191978423} & \num[round-precision=1]{916080431} & \num[round-precision=1]{3999010197} & \num[round-precision=1]{16712747961} & \num[round-precision=1]{68749088135} \\
        \bottomrule
      \end{tabular}
    \end{minipage}
  \end{subcaptionblock}
\end{table}
As a first test case, we consider a model problem on the space-time domain $\Omega\times I = [0,1]^2\times [0, 1]$ with prescribed solution given for the velcity $\boldsymbol v \colon \Omega\times I \to \mathbb{R}^2$ and pressure
$p \colon \Omega\times I \to \mathbb{R}$ by
\begin{subequations}\label{eq:conv-test}
\begin{align}
  \label{eq:conv-test-v}
  \boldsymbol v(\mathbf{x},\,t) &= \sin(t) \begin{pmatrix} \sin^2(\pi x) \sin(\pi y) \cos(\pi y) \\
                                        \sin(\pi x) \cos(\pi x) \sin^2(\pi y) \end{pmatrix},\\
\label{eq:conv-test-p} p(\mathbf{x},\,t) &= \sin(t) \sin(\pi x) \cos(\pi x)
  \sin(\pi y) \cos(\pi y)\,.
\end{align}
\end{subequations}
We set the kinematic viscosity to $\nu = 0.1$ and choose
the external force $\mathbf{f}$ such that the solution~\eqref{eq:conv-test}
satisfies~\eqref{Eq:SE}. The initial velocity is
prescribed as zero and homogeneous Dirichlet boundary conditions are imposed on
$\partial\Omega$ for all times
\[
  \boldsymbol v = \mathbf{0}\text{ on }\Omega\times \{0\},\quad
  \boldsymbol v = \mathbf{0}, \text{ on } \partial \Omega\times (0, T]\,.
\]
The space-time mesh $\mathcal{T}_{h}\otimes\mathcal{M}_{\tau}$ is a uniform
triangulation of the space-time domain $\Omega\times I$. We use discretizations
with varying polynomial degrees $r \in \{3,\, 4,\dots,\,8\}$ in space and
$k=r$ in time to test the convergence.

Table~\ref{tab:conv-stokes} shows the findings of our convergence study for
$r\in\{4,5\}$. The expected orders of convergence match with the experimental
orders. For a full account of all tests, we refer to \Cref{fig:conv-stokes}
in Appendix~\ref{Sec:CP}. The $L^2(0,T;L^2(\Omega)^d)$ velocity error
  does not always reach the ideal order $r+2$ due to the temporal polynomial
  degree $r$. However, the $L^2(0,T;H^1(\Omega)^d)$ norm exhibits the optimal
  order $r+1$, which justifies our choice of polynomial order $r$ in time.
Table~\ref{tab:iter-stokes} shows the number of GMRES iterations required for
convergence for these experiments. We compare the $hp$ STMG with a pure spatial
$h$-multigrid method and cell patches with VSPs. The $hp$ STMG method
exhibits superior robustness and significantly reduces the number of iterations,
as the polynomial degree and mesh refinement increase. In specific settings on
coarse meshes and low polynomial orders, the $hp$ STMG can be outperformed by
other strategies (e.\,g.\ pure $h$-multigrid, $h$-multigrid in space and
$p$-multigrid in time only). However, it consistently outperforms them on finer
meshes and higher polynomial degrees. In terms of total elements in the
  smoother per subinterval $I_n$, denoted by $N_{\text{sm}}^{\text{el}}$, the savings of $hp$ STMG may appear moderate. The main advantage
  lies in the smaller block sizes from the polynomial multigrid and the reduced
  GMRES iterations, which lower the number of smoother executions.
Further, although polynomial coarsening $r\rightarrow r-1$ can reduce GMRES
iterations and improve $p$-robustness, it yields no significant gains in
wall-clock time due to the slower decrease of the block size. The reduction of
the block size through polynomial coarsening increases the efficiency
significantly. In particular, halving the polynomial degree
$r\rightarrow \tfrac{r}{2}$ provides the best improvements in solver
performance. 
  The VSP smoother shows no advantage in the $hp$-STMG setting, as the
  velocity-pressure coupling is already captured within cell-based patches.
  Moreover, it is costly, with the total number of elements in a patch often
  exceeding that of cell-based patches by an order of magnitude, while providing
  only marginal performance gains. Consequently, we omit it in large-scale
  simulations. For pairs with a wider pressure-velocity stencil, such as
  Taylor-Hood, VSPs can offer benefits.

We note that only a single smoothing step, i.\,e.\ \(\nu_1=\nu_2=1\) is
performed on all levels. While additional smoothing steps could reduce the
number of GMRES iterations and improve the $h$- and $p$-robustness, it may not
improve the time to solution for matrix-free methods. We use a matrix-based
smoother~\eqref{vanka}, so keeping the number of smoothing steps small and
reducing the complexity is the most important part of the overall performance,
see~\cite{margenbergSpaceTimeMultigridMethod2024a}. In the present section, we
achieve excellent $h$-robustness, but not full $p$-robustness. We revisit this
topic in the following section and address whether an increase in smoothing
steps improves the $p$-robustness.

\subsection{Lid-driven cavity flow}
\begin{table}[htb]
  \centering
  \caption{Average number of GMRES iterations per subproblem
    $\overline{n}_{\text{iter}}$ for the number of smoothing steps
    $\nu_1=\nu_2=n_{\text{sm}}$, polynomial degrees $r$ and refinements $c$. We include
    the number of global space-time cells (\#
    st-cells) and smoother elements $N_{\text{sm}}^{\text{el}}$.}\label{tab:lid-iter-stokes}
  \sisetup{round-precision=3}
  \footnotesize
  \begin{tabular}{ll|ll|rr|rr|rr}
    \toprule
    &                &\multicolumn{2}{c|}{$N_{\text{sm}}^{\text{el}}$}& \multicolumn{2}{c|}{$n_{\text{sm}}=1$}                              &\multicolumn{2}{c|}{$n_{\text{sm}}=2$}                  &\multicolumn{2}{c}{$n_{\text{sm}}=4$}\\
    $c$&\# st-cells   & {$r=2$}   & {$r=3$} & {$r=2$}  & {$r=3$}                 & {$r=2$}   & {$r=3$}   & {$r=2$}   & {$r=3$}\\
  \midrule
$4$&\num{1048576}   &\num[round-precision=2]{8.03619e+08} &\num[round-precision=2]{6.05557e+09}&18.14  &28.09  & 10.56  & 16.22  & 6.83  & 10.66 \\
$5$&\num{16777216}  &\num[round-precision=2]{6.43034e+09} &\num[round-precision=2]{4.84482e+10}&16.97  &25.72  &  9.32  & 13.79  & 5.83  &  9.17 \\
$6$&\num{268435456} &\num[round-precision=2]{5.14441e+10} &\num[round-precision=2]{3.8759e+11}&14.82  &21.95  &  7.60  & 11.26  & 4.86  &  7.56 \\
$7$&\num{4294967296}&\num[round-precision=2]{4.11554e+11} &\num[round-precision=2]{3.10072e+12}&12.52  &18.41  &  6.39  &  9.25  & 3.88  &  6.27 \\
  \bottomrule
  \end{tabular}
  \caption{Throughput $\theta$~\eqref{eq:throughput} for different values of $n_{\text{sm}}$, $r$ and $c$.}\label{tab:lid-throughput-stokes}
  \setlength{\tabcolsep}{5pt}
  \begin{tabular}{ll|rr|rr|rr}
    \toprule
    &                & \multicolumn{2}{c|}{$n_{\text{sm}}=1$}                              &\multicolumn{2}{c|}{$n_{\text{sm}}=2$}                  &\multicolumn{2}{c}{$n_{\text{sm}}=4$}\\
    $c$&\# st-cells   & {$r=2$}  & {$r=3$}                 & {$r=2$}   & {$r=3$}   & {$r=2$}   & {$r=3$}\\
    \midrule
  $4$&\num{1048576}    &\eval{302520576/149.4}    &\eval{927534080/629.4}  &\eval{302520576/158.4}  &\eval{927534080/699.2} &\eval{302520576/185} &\eval{927534080/883.8}\\
  $5$&\num{16777216}   &\eval{4708913664/329}     &\eval{14531434496/1246}  &\eval{4708913664/331.2}  &\eval{14531434496/1285} &\eval{4708913664/375.5}  &\eval{14531434496/1641}\\
  $6$&\num{268435456}  &\eval{74307412992/958.8}  &\eval{230058635264/5066}  &\eval{74307412992/917.3}  &\eval{230058635264/6012}  &\eval{74307412992/1105} &\eval{230058635264/6647}\\
  $7$&\num{4294967296} &\eval{1180701050880/6502} &\eval{3661497393152/45990} &\eval{1180701050880/7483}  &\eval{3661497393152/47720} &\eval{1180701050880/7977} &\eval{3661497393152/63030}\\
  \bottomrule
  \end{tabular}
\end{table}
\begin{figure}[htbp]
\includegraphics{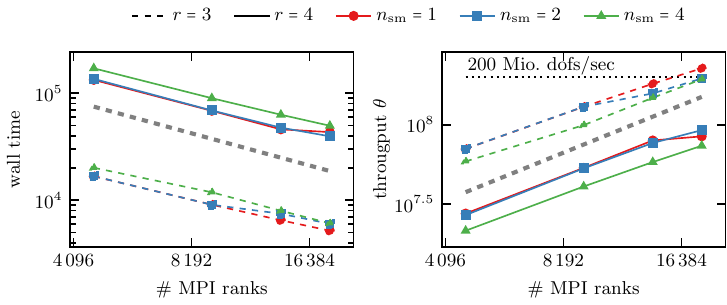}
\caption{\label{fig:lid-strong-scale}Strong scaling test results for the STMG
  algorithm with varying numbers of smoothing steps. The left plot shows the
  time to solution over the number of MPI processes. The dashed gray lines
  indicate the optimal scaling. The right plot depicts the degrees of freedom
  (dofs) processed per second over the number of MPI processes.}
\end{figure}
\begin{figure}[htbp]
  \includegraphics{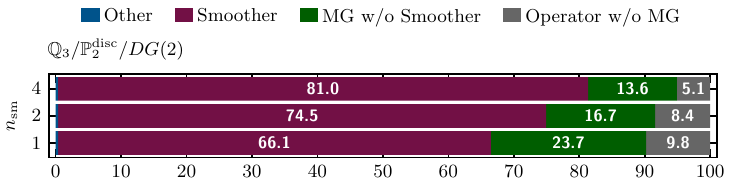}
  \includegraphics{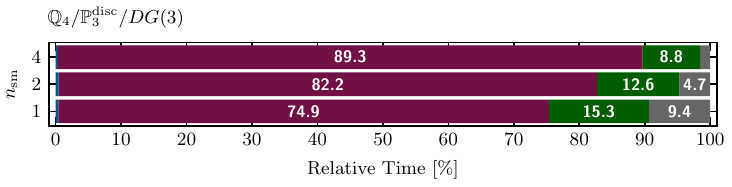}
  \caption{\label{fig:lid-strong-scale-rel}Time spent in different parts of the
    lid-driven cavity flow simulation, executed on
    $\num[scientific-notation=false,round-precision=0]{18432}$ MPI processes.
    The simulations were conducted for $c=7$, $r\in\{3,\,4\}$,
    $n_{\text{sm}}\in\{1,\,2,\,4\}$.}
\end{figure}
We now study the more sophisticed benchmark problem of lid-driven cavity flow.
The space-time mesh $\mathcal{T}_{h}\times\mathcal{M}_{\tau}$ is a uniform triangulation of the space-time domain $\Omega\times I=[0,\,1]^3\times [0,\,8]$, refined globally $c$ times. A Dirichlet profile
\(\mathbf v_{D}\) is prescribed at the upper boundary
\(\Gamma_{D}=[0,\,1]^{2}\times \{1\}\subset \partial \Omega\) as
\begin{equation}
    \mathbf v_{D}(x,\,y,\,z,\,t)=\sin\left(\tfrac{\pi}{4} t\right) \text{  on  } \Gamma_{D}\times [0,\,8]\,.
\end{equation}
On the other boundaries, denoted by
\(\Gamma_{\text{wall}}=\partial \Omega \setminus\Gamma_{D}\) we use no-slip
boundary conditions. We employ discretizations with different polynomial degrees
$r \in \{2,\, 3\}$ in space. For the time discretization we set $k=r$. For the
strong scaling test shown in \Cref{fig:lid-strong-scale} we set $c=7$, which
results in $\num[scientific-notation=false,round-precision=0]{2048}$ time cells
and $\num{2097152}$ space cells.
This
configuration yields \num[scientific-notation=true]{446960131} and \num[scientific-notation=true]{192171395} spatial degrees
of freedom for \(r=3\) and \(r=2\), respectively, and
\num[scientific-notation=true]{3661497393152} and \num[scientific-notation=true]{1180701050880} global space-time dofs. The
resulting local linear systems each involve \num[scientific-notation=true]{1787840524} unknowns for
\(r=3\) and \num[scientific-notation=true]{576514185} unknowns for \(r=2\). The average number of GMRES iterations
$\overline{n}_{\text{iter}}$ for different values of
$\nu_1=\nu_2\eqcolon n_{\text{sm}}\in \{1,\,2,\,4\}$ (cf.~\eqref{Alg:VGMGPre},~\eqref{Alg:VGMGPost}), $r \in \{2,\, 3\}$ and
$c\in \{4,\,5,\,6,\,7\}$ are collected in Table~\ref{tab:lid-iter-stokes}. To
compare the computational efficiency, we measure the
\emph{throughput}. Let \( W_{\text{total}}(n_{\text{sm}}, c, r) \) be the total walltime,
and \( N_{\text{dof}}(c, r) \) the total number of degrees of freedom. The
throughput $\theta(n_{\text{sm}}, c, r)$ is defined as
\begin{equation}
\label{eq:throughput}
\theta(n_{\text{sm}}, c, r) = \frac{N_{\text{dof}}(c, r)}{W_{\text{total}}(n_{\text{sm}}, c, r)}.
\end{equation}
In addition to the scaling test and performance of the iterative solver, we
verify the convergence of our discretization by studying a goal quantity in
Appendix~\ref{sec:pdiff}. In Table~\ref{tab:lid-throughput-stokes}, we summarize
\( \theta \) using walltimes obtained on
$\num[scientific-notation=false,round-precision=0]{13824}$ MPI ranks in the same
configurations as presented in Table~\ref{tab:lid-iter-stokes}. Increasing the
number of smoothing steps, \(n_{\text{sm}}\), reduces the GMRES iterations, but
the higher cost per iteration leads to increased wall times. A larger
\(n_{\text{sm}}\) also improves \(p\)-robustness, which is satisfactory but
could be improved. Overall, the Vanka smoother~\eqref{vanka} is effective.
However, cell-wise direct solves introduce significant overhead, especially at
higher polynomial degrees. In \Cref{fig:lid-strong-scale-rel} we show the
relative execution time spent in different parts of the program within the
strong scaling test for different $n_{\text{sm}}$. The \textit{MG w/o Smoother}
segment corresponds to the $hp$ STMG without its smoothing steps, i.e.\ operator
evaluations and grid transfers. Its contribution decreases as the number of
smoothing steps increases. In contrast, the \textit{Smoother} segment,
consistently dominates the wall time. Its cost depends on the polynomial order
and continues to increase with $n_{\text{sm}}$. The \textit{Operator w/o MG}
part covers operator evaluations performed outside the $hp$ STMG preconditioner,
and the \textit{Other} segment includes the time spent on source term assembly,
goal quantity evaluation and tasks between time steps. The absolute time of
these segments remain constant, resulting in a decrease in relative time with
increasing smoothing steps.

\section{Conclusions}
\label{sec:conclusions}
We present an $hp$ multigrid approach for tensor-product space-time finite
element discretizing the Stokes equations. The method exhibits optimal
$h$-robustness and satisfactory $p$-robustness. Even when embedded in a
time-marching scheme, the proposed $hp$ space-time multigrid method achieves
notable efficiency gains through combined space-time coarsening. The smoother is
effective but expensive and represents the main computational bottleneck for
higher order discretizations. Replacing the direct solver within the smoother
with a more efficient local method becomes crucial to avoid the escalating cost
at large \(p\). The direct solves raise concerns about memory usage. Iterative
local solvers might be less resource intensive. Block-diagonal or approximate
factorization approaches, e.\,g.\ a diagonal Vanka~\cite{john_numerical_2000}, might preserve efficiency without incurring the cost of
a direct solver. The extension of temporal decoupling~\cite{munchStageParallelFullyImplicit2023} to the local systems in the
smoother could further reduce the costs. Despite the current
limitations, the method performs well and achieves throughput over 200 millions
of degrees of freedom per second on problems with trillions of global degrees of
freedom. It outperforms existing matrix-based implementations by orders of
magnitude (cf.~\cite{anselmannGeometricMultigridMethod2023}). The proposed
matrix-free $hp$ multigrid method for tensor-product space-time finite element
discretizations is highly efficient and scalable, making it a promising candidate for
large-scale problems in fluid mechanics, fluid-structure interaction, and
dynamic poroelasticity. An extension to these problems, particularly nonlinear
ones (e.g.\ Navier--Stokes), is part of future work.

\appendix
\setcounter{table}{0}
\renewcommand{\thetable}{\Alph{section}.\arabic{table}}
\setcounter{figure}{0}
\renewcommand{\thefigure}{\Alph{section}.\arabic{figure}}

\section{$h$- and $p$-Multigrid Transfer Operators for Space-Time}
\label{App:TOMG}
Here, we explicitly present the construction and definition of the restriction and prolongation operators introduced in \Cref{Sec:GTO} for our $h$- and $p$-multigrid concept. We address the transfer operators for the refinement and coarsening of the space and time mesh of the algebraic tensor product
spaces~\eqref{Eq:STSMG}, referred to as $h$-multigrid. We start with the spatial mesh. For the nested finite element spaces~\eqref{Eq:NFES} and~\eqref{Eq:NFESp+}
with bases according to~\eqref{Eq:BasVhQh} we define the isomorphisms which map
the degrees of freedom to the finite element spaces as (cf.~\cite{olshanskiiMultigridAnalysis2012})
\begin{subequations}
	\label{Def:IsoAlg}
	\begin{alignat}{7}
		\boldsymbol{\mathcal R}^{r+1}_s : & \;  \mathbb R^{M^{\boldsymbol v}_{r+1,s}} & \to & \, \boldsymbol V_s^{r+1}\,, & \qquad && \boldsymbol{\mathcal R}^{r+1}_s \boldsymbol U_s & = \sum_{i=0}^{M^{\boldsymbol v}_{r+1,s}} U_s^i \boldsymbol  \chi_{i,s}^{\boldsymbol{v}}\,,\\
		\mathcal S^{r}_s : &  \; \mathbb R^{M_{r,s}^p} & \to & \; Q_s^{r,+}\,, && \qquad & \mathcal S^{r}_s \boldsymbol  P_s  & = \sum_{i=0}^{M_{r,s}^p} P_s^i \chi^p_{i,s}\,.
	\end{alignat}
\end{subequations}
We note that $Q_s^{r} = \{\mathcal S^{r}_s \boldsymbol P_s \mid \boldsymbol P_s \in (\boldsymbol M_{h,s}^p \boldsymbol 1)^\perp \}\}$ with $\boldsymbol 1=(1,\ldots,1)^\top\in \R^{M_{r,s}^p}$ and pressure mass matrix $\boldsymbol M_{r,s}^p\in \R^{M_{r,s}^p,M_{r,s}^p}$  of~\eqref{Eq:DefMAB_2}. The orthogonality condition $\langle p,1\rangle_{L^2(\Omega)}=0$ for $p\in Q_h^{r,s}$ corresponds to the orthogonality $\langle \boldsymbol P,\boldsymbol M_{h,s}^p\boldsymbol 1\rangle_{\mathbb R^{M_{r,s}^p}} =0$ for $\boldsymbol P \in \mathbb R^{M_{r,s}^p}$. For $s=0,\ldots,S$ and $r\in \N$, we let
\begin{equation}
	\label{Eq:DefRsSs}
	\boldsymbol R^{r+1}_s \coloneq \mathbb R^{M^{\boldsymbol v}_{r+1,s}}  \quad \text{and} \quad \boldsymbol S^{r}_s \coloneq (\boldsymbol M_{h,s}^p \boldsymbol 1)^\perp\,.
\end{equation}
For prolongation and restriction in space and with~\eqref{Def:IsoAlg}, we use the canonical choice
\begin{subequations}
	\label{Eq:DefPrRe}
	\begin{alignat}{5}
		\label{Eq:DefTps}
		&\begin{aligned}
			& \boldsymbol T^{r}_{s-1,s}:  \boldsymbol R^{r+1}_{s-1} \times \boldsymbol S^{r}_{s-1}	 \to  \boldsymbol R^{r+1}_s \times \boldsymbol S^{r}_s\,, \\[3pt]
			& \boldsymbol T^{r}_{s-1,s} = (\boldsymbol{\mathcal R}^{r+1}_s)^{-1}\circ \boldsymbol{\mathcal R}^{r+1}_{s-1} \times (\mathcal S^{r}_s)^{-1}\circ\mathcal S^{r}_{s-1}\,,
		\end{aligned}
		\\[3pt]
		\label{Eq:DefTrs}
		& \begin{aligned}
			&\boldsymbol T^{r}_{s,s-1}:  \boldsymbol R^{r+1}_{s} \times \boldsymbol  S^{r}_{s}	\to  \boldsymbol R^{r+1}_{s-1} \times \boldsymbol S^{r}_{s-1}\,,\\[3pt]
			& \boldsymbol T^{r}_{s,s-1}  = (\boldsymbol{\mathcal R}^{r+1}_{s-1})^\ast\circ ((\boldsymbol{\mathcal R}^{r+1}_{s})^\ast)^{-1} \times (\mathcal S^r_{s-1})^{\ast}\circ ((\mathcal S^r_{s})^\ast)^{-1}\,.
		\end{aligned}
	\end{alignat}
\end{subequations}
Both operators, $\boldsymbol T^r_{s-1,s}$ and $\boldsymbol T^r_{s,s-1}$, keep the pressure in the correct subspace. For $\boldsymbol{Z}=(\boldsymbol  V,\boldsymbol  P)^\top\in (\boldsymbol R_{s-1}^{r+1})^{k+1} \times (\boldsymbol S_{s-1}^r)^{k+1}$, with $\boldsymbol V \in  (\boldsymbol R_{s-1}^{r+1})^{k+1}$ and $\boldsymbol P \in  (\boldsymbol S_{s-1}^r)^{k+1}$, we define the prolongation $\boldsymbol{\widetilde T}^r_{s-1,s}:  (\boldsymbol R_{s-1}^{r+1})^{k+1} \times (\boldsymbol S_{s-1}^r)^{k+1} \to  (\boldsymbol R_{s}^{r+1})^{k+1} \times (\boldsymbol S_{s}^r)^{k+1}$ for the product spaces $(\boldsymbol R_{s-1}^{r+1})^{k+1}$ and $(\boldsymbol S_{s-1}^r)^{k+1}$ by componentwise application of~\eqref{Eq:DefTps},
\begin{equation}
	\label{Eq:DeftTps}
	\begin{aligned}
		\boldsymbol{\widetilde T}_{s-1,s}^r \boldsymbol{Z} & \coloneq (\boldsymbol T_{s-1,s}^{r,\boldsymbol v} \boldsymbol V^1, \ldots,\boldsymbol T_{s-1,s}^{r,\boldsymbol v}\boldsymbol V^{k+1},  \boldsymbol T_{s-1,s}^{r,p} \boldsymbol P^1, \ldots,  \boldsymbol T_{s-1,s}^{r,p} \boldsymbol P^{k+1})^\top\\[3pt]
		& =  \begin{pmatrix} \boldsymbol E_{k+1}  \otimes \boldsymbol T_{s-s,s}^{r,\boldsymbol v} & \boldsymbol 0 \\  \boldsymbol 0 & \boldsymbol E_{k+1}  \otimes \boldsymbol T_{s-1,s}^{r,p}  \end{pmatrix}  \begin{pmatrix} \boldsymbol V\\ \boldsymbol P \end{pmatrix}\,,
	\end{aligned}
\end{equation}
where $ \boldsymbol T_{s-1,s}^{r,\boldsymbol v}:\boldsymbol R_{s-1}^{r+1}  \to
\boldsymbol R_{s}^{r+1}$ and $ \boldsymbol T_{s-1,s}^{r,p}:\boldsymbol S_{s-1}^r
\to \boldsymbol S_{s}^r $ are the velocity and pressure parts of the
prolongation $ \boldsymbol T_{s,s-1}^r$ defined in~\eqref{Eq:DefTps}. The corresponding restriction operator $\boldsymbol{\widetilde T}^r_{s,s-1}: (\boldsymbol R_{s}^{r+1})^{k+1} \times (\boldsymbol S_{s}^{r})^{k+1} \to  (\boldsymbol R_{s-1}^{r+1})^{k+1} \times (\boldsymbol S_{s-1}^{r})^{k+1}$  is defined analogously by
\begin{equation}
	\label{Eq:DeftTrs}
	\begin{aligned}
		\boldsymbol{\widetilde T}_{s,s-1}^r \boldsymbol{Z} & \coloneq (\boldsymbol T_{s,s-1}^{r,\boldsymbol v} \boldsymbol V^1, \ldots,\boldsymbol T_{s,s-1}^{r,\boldsymbol v}\boldsymbol V^{k+1},  \boldsymbol T_{s,s-1}^{r,p} \boldsymbol P^1, \ldots,  \boldsymbol T_{s,s-1}^{r,p} \boldsymbol P^{k+1})^\top\\[3pt]
		& = \begin{pmatrix}  \boldsymbol E_{k+1}  \otimes  \boldsymbol T_{s,s-1}^{r,\boldsymbol v} & \boldsymbol 0 \\  \boldsymbol 0 &  \boldsymbol E_{k+1}  \otimes \boldsymbol T_{s,s-1}^{r,p}  \end{pmatrix}  \begin{pmatrix} \boldsymbol V\\ \boldsymbol P \end{pmatrix}\,,
	\end{aligned}
\end{equation}
where
$ \boldsymbol T_{s,s-1}^{r,\boldsymbol v}:\boldsymbol R_{s}^{r+1} \to
\boldsymbol R_{s-1}^{r+1}$ and
$ \boldsymbol T_s^{r,p}:\boldsymbol S_{s}^r \to \boldsymbol S_{s-1}^r $ are the
velocity and pressure parts of the restriction $ \boldsymbol T_{s,s-1}^r$
defined in~\eqref{Eq:DefTrs}. For global vectors
$\boldsymbol{Z}=(\boldsymbol Z_1\ldots,\boldsymbol Z_{N_l})^\top\in
\big((\boldsymbol R_{s-1}^{r+1})^{k+1} \times (\boldsymbol
S_{s-1}^r)^{k+1}\big)^{N_l}$ with
$\boldsymbol Z_n = (\boldsymbol V_n,\boldsymbol P_n)^\top\in (\boldsymbol
R_{s-1}^{r+1})^{k+1} \times (\boldsymbol S_{s-1}^r)^{k+1}$, for
$n=1,\ldots,N_l$, and with
$\boldsymbol 1_{N_l} \coloneq (1,\ldots,1)\in \mathbb R^{N_l}$ according to~\eqref{Eq:DefSubvec_1} to~\eqref{Eq:DefX}, prolongation
and restriction are then defined by
\begin{equation}
	\label{Eq:DefgT}
	\boldsymbol{\overline T}_{s-1,s}^r \boldsymbol{Z} \coloneq  \big((\boldsymbol{1}_{N_l} \otimes \boldsymbol{\widetilde T}_{s-1,s}^r) \boldsymbol{Z}\big)^\top \quad \text{and} \quad
	\boldsymbol{\overline T}_{s,s-1}^r \boldsymbol{Z} \coloneq  \big((\boldsymbol{1}_{N_l} \otimes \boldsymbol{\widetilde T}_{s,s-1}^r) \boldsymbol{Z}\big)^\top\,.
\end{equation}
Next, we introduce the grid transfer operator for the hierarchy of temporal
meshes. Similarly to~\eqref{Def:IsoAlg}, for temporal mesh level $l=0,\ldots,L$ we define the isomorphism, on spatial mesh level $s=0,\ldots ,S$, with $D^{k,r}_{l,s} \coloneq N_l (k+1)(M_{r+1,s}^{\boldsymbol v}+M_{r,s}^p)$ as
\begin{equation}
	\label{Def:IsoTme}
	\boldsymbol{\mathcal  I}_l^k  : \mathbb R^{D^{k,r}_{l,s}}\rightarrow Y_l^{k}(I) \otimes \mathbb R^{{M_{r+1,s}^{\boldsymbol v}+M_{r,s}^p} }\,,  \qquad  \boldsymbol{\mathcal  I}_l^k  \boldsymbol W = \sum_{n=1}^{N_l} \sum_{a=1}^{k+1} \begin{pmatrix} \boldsymbol W^{l,\boldsymbol v}_{n,a} \\[3pt] \boldsymbol W^{l,p}_{n,a} \end{pmatrix}\varphi_{n,a}^l\,.
\end{equation}
Here, we substructured $\boldsymbol W$ according to
\begin{equation}
	\label{Eq:SubStrGloVec}
	\begin{aligned}
		\boldsymbol W = \big(& \boldsymbol W^{l,\boldsymbol v}_{1,1},\ldots, \boldsymbol W^{l,\boldsymbol v}_{1,k+1},\boldsymbol W^{l,p}_{1,1},\ldots, \boldsymbol W^{l,p}_{1,k+1},\ldots, \\[3pt]
		& \boldsymbol W^{l,\boldsymbol v}_{N_l,1},\ldots ,\boldsymbol W^{l,\boldsymbol v}_{N_l,k+1} , \boldsymbol W^{l,p}_{N_l,1},\ldots ,\boldsymbol W^{l,p}_{N_l,k+1}\big)^\top\,,
	\end{aligned}
\end{equation}
with $\boldsymbol W_{n,a}^{l,\boldsymbol v}\in \mathbb R^{M_{r+1,s}^{\boldsymbol v}}$ and $\boldsymbol W_{n,a}^{l,p}\in \mathbb R^{M_{r,s}^{p}}$ for $n=1,\ldots,N_l$ and $a=1,\ldots,k+1$. For prolongation and restriction we use the canonical choices again,
\begin{subequations}
	\label{Eq:DeftPrRe}
	\begin{alignat}{5}
		\boldsymbol I^k_{l-1,l}: & \; \mathbb R^{D^{k,r}_{l-1,s}} \to \mathbb R^{D^{k,r}_{l,s}} \,, & \qquad \boldsymbol I^{k}_{l-1,1} & =  (\boldsymbol{\mathcal  I}^k_{l})^{-1} \circ \boldsymbol{\mathcal  I}^k_{l-1}\,,\\[3pt]
		\boldsymbol I_{l,l-1}^k: &  \; \mathbb R^{D^{k,r}_{l,s}}  \to  \mathbb R^{D^{k,r}_{l-1,s}} \,, & \qquad \boldsymbol I_{l,l-1}^k & = (\boldsymbol{\mathcal  I}^k_{l-1})^{\ast}  \circ	((\boldsymbol{\mathcal  I}^k_{l})^{\ast})^{-1}\,.
	\end{alignat}
\end{subequations}
Space-time prolongation $\boldsymbol T_{l-1,l;s-1,s}^{k,r}: \mathbb R^{D^{k,r}_{l-1,s-1}}\to \mathbb R^{D^{k,r}_{l,s}}$ and restriction $\boldsymbol T^{k,r}_{l,l-1;s,s-1}: \mathbb R^{D^{k,r}_{l,s}}\to \mathbb R^{D^{k,r}_{l-1,s-1}}$  are then defined by the concatenation of $\boldsymbol I_{l-1,l}^k$ and $\boldsymbol{\overline T}^r_{s-1,s}$ and of $\boldsymbol I_{l,l-1}^r$ and $\boldsymbol{\overline T}^r_{s,s-1}$, respectively, as
\begin{equation}
	\label{Eq:DefTlspr}
	\boldsymbol T_{l-1,l;s-1,s}^{k,r} \coloneq \boldsymbol I_{l-1,l}^k \circ  \boldsymbol{\overline T}^r_{s-1,s} \qquad \text{and} \qquad
	\boldsymbol T_{l,l-1;s,s-1}^{k,r}  \coloneq \boldsymbol I_{l,l-1}^k \circ \boldsymbol{\overline T}^r_{s,s-1}\,.
\end{equation}
\begin{remark}
  If the time-marching process~\eqref{Eq:AlgSysTMS} is applied, we use
  $hp$-multigrid in space and $p$-multigrid in time and omit
  $h$-multigrid in time.
  Then, the operators $\boldsymbol{I}_{l-1,l}^k$ and
  $\boldsymbol{I}_{l,l-1}^k$ in~\eqref{Eq:DefTlspr} reduce to the identity. The
  $hp$-STMG (\Cref{alg:stmg}) remains well-defined. This also holds for macro
  time steps $1<\widetilde{N} < N$ (cf. Remark~\ref{rem:GLSP}).
\end{remark}

We define the transfer operators for coarsening and prolongation of
the polynomial degrees $k,r \in \mathbb{N}$ in~\eqref{Eq:STSMG}, referred to as
$p$-multigrid. For simplicity, we assume
\begin{equation}
	\label{Assmp:MG2}
	k = 2^K \quad \text{for some } K \in \mathbb{N},\quad r = 2^R \quad \text{for some } R \in \mathbb{N}\,,
\end{equation}
which facilitates bisection coarsening $k \mapsto \lfloor k/2 \rfloor$ and
doubling for prolongation. Since the floor function is not invertible, one would
otherwise need an accounting vector for the polynomial orders. This notational
overhead is avoided by~\eqref{Assmp:MG2}. The general case of arbitrary
$k,r \in \mathbb{N}$ is treated in Algorithm~\ref{alg:ConstructHierarchy}. The
bisection strategy is motivated by computational studies showing that it strikes
a favorable balance between two-level ($k \mapsto 1$) and decrement-by-one
($k \mapsto k-1$) coarsening;
see~\cite[Section~3.2.2]{fehnHybridMultigridMethods2020} for a review.

For prolongation of the spatial polynomial order $r$, i.e.\
$\frac{r}{2}\mapsto r$ in~\eqref{Eq:NFES}, and restriction, i.e.\
$r\mapsto\frac{r}{2}$, we let the transfer operators (for $R\geq 2$)
\begin{equation}
	\label{Eq:DefPrRer}
	\boldsymbol T_s^{\frac{r}{2},r} :  \; \boldsymbol R_{s}^{\frac{r}{2}}  \times  \boldsymbol S_s^{\frac{r}{2}-1}   \to    \boldsymbol R_s^{r} \times \boldsymbol S_s^{r-1}\,, \qquad
	\boldsymbol T_s^{r,\frac{r}{2}} :  \; \boldsymbol R_s^{r}  \times  \boldsymbol S_s^{r-1}   \to    \boldsymbol S_s^{\frac{r}{2}} \times \boldsymbol S_s^{\frac{r}{2}-1}
\end{equation}
for the spaces~\eqref{Eq:DefRsSs} be defined analogously to~\eqref{Eq:DefPrRe}
along with~\eqref{Def:IsoAlg}. The prolongation $\boldsymbol{\overline
	T}_s^{\frac{r}{2},r}$ and restriction $\boldsymbol{\overline T}_s^{r,\frac{r}{2}}$ extend
$\boldsymbol{T}_s^{\frac{r}{2},r}$  and $\boldsymbol{T}_s^{r,\frac{r}{2}}$ to the global vector
of unknowns~\eqref{Eq:SubStrGloVec} along the lines of~\eqref{Eq:DeftTps}
and~\eqref{Eq:DeftTrs}.
Finally, we define the prolongation and restriction of the temporal polynomial
order \( k\) in~\eqref{Eq:DefYkl}. For
prolongation and restriction in time of the $p$-multigrid method,
\begin{equation}
	\label{Eq:DefPrRek}
	\boldsymbol I_l^{\frac{k}{2},k} :  \; \mathbb R^{D^{\frac{k}{2},r}_{l,s}} \to \mathbb R^{D^{k,r}_{l,s}} \,, \qquad
	\boldsymbol I_l^{k,\frac{k}{2}} :  \; \mathbb R^{D^{k,r}_{l,s}}  \to  \mathbb R^{D^{\frac{k}{2},r}_{l,s}}
\end{equation}
are defined analogously to~\eqref{Eq:DeftPrRe} along with~\eqref{Def:IsoTme}. Space-time combined grid transfer operations of the $p$-multigrid method are then constructured by concatenation of the operators in~\eqref{Eq:DefPrRer} and~\eqref{Eq:DefPrRek}, similarly to~\eqref{Eq:DefTlspr}, such that
\begin{equation}
	\label{Eq:DefoTsl}
	\boldsymbol T_{l,s}^{\frac{k}{2},k;\frac{r}{2}r,r} \coloneq \boldsymbol I_{l}^{ \frac{k}{2},k} \circ \boldsymbol{\overline T}^{\frac{r}{2}r,r}_{s} \qquad \text{and} \qquad
	\boldsymbol T_{l,s}^{k,\frac{k}{2};r,\frac{r}{2}}  \coloneq \boldsymbol I_{l}^{k,\frac{k}{2}} \circ \boldsymbol{\overline T}^{r,\frac{r}{2}}_{s}\,.
\end{equation}

\section{Subroutines in the Multigrid algorithm}
The subroutines of \Cref{alg:stmg} are summarized now.	
\begin{algorithm}[H]
	\caption{\textsc{P-Residual-Restrict}}
	\label{alg:pmg}
  \small
	\textbf{Input:} fine-level iterate \( \boldsymbol{X}^{k,r;\nu_1}_{l,s} \), right-hand side \( \boldsymbol{B}^{k,r}_{l,s} \), system matrix \( \boldsymbol{S}^{k,r}_{l,s} \) \\
	\begin{enumerate}\itemsep1pt \parskip0pt \parsep0pt
		\item[\textbf{1.}] Compute the fine-level residual:
		\(\displaystyle
		\boldsymbol{R}^{k,r}_{l,s} = \boldsymbol{B}^{k,r}_{l,s} - \boldsymbol{S}^{k,r}_{l,s} \boldsymbol{X}^{k,r;\nu_1}_{l,s}.
		\)
		
		\item[\textbf{2.}] \textbf{If} \( r > K \): restrict and solve
		\(\displaystyle
		\boldsymbol{B}^{k,\frac{r}{2}}_{l,s} = \boldsymbol{T}^{k,k;r,\frac{r}{2}}_{l,s} \boldsymbol{R}^{k,r}_{l,s},
		\)
		\[
		\text{Find } \boldsymbol{Y}_{{p}}^{k,\frac{r}{2}} \in \mathbb{R}^{D^{k,\frac{r}{2}}_{l,s}} \text{ such that } 
		\boldsymbol{S}^{k,\frac{r}{2}}_{l,s} \boldsymbol{Y}_{{p}}^{k,\frac{r}{2}} = \boldsymbol{B}^{k,\frac{r}{2}}_{l,s}.
		\]
		
		\item[\textbf{3.}] \textbf{Else}, restrict and solve
		\(\displaystyle
		\boldsymbol{B}^{\frac{k}{2},\frac{r}{2}}_{l,s} = \boldsymbol{T}^{k,\frac{k}{2};r,\frac{r}{2}}_{l,s} \boldsymbol{R}^{k,r}_{l,s},
		\)
		\[
		\text{Find } \boldsymbol{Y}_{{p}}^{\frac{k}{2},\frac{r}{2}} \in \mathbb{R}^{D^{\frac{k}{2},\frac{r}{2}}_{l,s}} \text{ such that } 
		\boldsymbol{S}^{\frac{k}{2},\frac{r}{2}}_{l,s} \boldsymbol{Y}_{{p}}^{\frac{k}{2},\frac{r}{2}} = \boldsymbol{B}^{\frac{k}{2},\frac{r}{2}}_{l,s}.
		\]
		\item[\textbf{4.}] \textbf{Return} \( \boldsymbol{Y}_{{p}}^{k',\frac{r}{2}} \)
	\end{enumerate}
\end{algorithm}
\vspace*{-1ex}
\begin{algorithm}[H]
	\caption{\textsc{H-Residual-Restrict}}
	\label{alg:hmg}
  \small
	\textbf{Input:} fine-level iterate \( \boldsymbol{X}^{1,1;\nu_1}_{l,s} \), right-hand side \( \boldsymbol{B}^{1,1}_{l,s} \), system matrix \( \boldsymbol{S}^{1,1}_{l,s} \) \\
	\begin{enumerate}\itemsep1pt \parskip0pt \parsep0pt
		\item[\textbf{1.}] Compute the fine-level residual:
		\(\displaystyle
		\boldsymbol{R}^{1,1}_{l,s} = \boldsymbol{B}^{1,1}_{l,s} - \boldsymbol{S}^{1,1}_{l,s} \boldsymbol{X}^{1,1;\nu_1}_{l,s}.
		\)
		
		\item[\textbf{2.}] \textbf{If} \( s > L \): restrict and solve
		\(\displaystyle
		\boldsymbol{B}^{1,1}_{l,s-1} = \boldsymbol{T}^{1,1}_{l,l;s,s-1} \boldsymbol{R}^{1,1}_{l,s},
		\)
		\[
		\text{Find } \boldsymbol{Y}^{{h}}_{l,s-1} \in \mathbb{R}^{D^{1,1}_{l,s-1}} \text{ such that } 
		\boldsymbol{S}^{1,1}_{l,s-1} \boldsymbol{Y}^{{h}}_{l,s-1} = \boldsymbol{B}^{1,1}_{l,s-1}.
		\]
		
		\item[\textbf{3.}] \textbf{Else}, restrict and solve
		\(\displaystyle
		\boldsymbol{B}^{1,1}_{l-1,s-1} = \boldsymbol{T}^{1,1}_{l,l-1;s,s-1} \boldsymbol{R}^{1,1}_{l,s},
		\)
		\[
		\text{Find } \boldsymbol{Y}^{{h}}_{l-1,s-1} \in \mathbb{R}^{D^{1,1}_{l-1,s-1}} \text{ such that } 
		\boldsymbol{S}^{1,1}_{l-1,s-1} \boldsymbol{Y}^{{h}}_{l-1,s-1} = \boldsymbol{B}^{1,1}_{l-1,s-1}.
		\]
		
		\item[\textbf{4.}] \textbf{Return} \( \boldsymbol{Y}^{{h}}_{l',s-1} \)
	\end{enumerate}
\end{algorithm}
\vspace*{-1ex}
\begin{algorithm}[H]
	\caption{\textsc{H-Prolongate-Correct}}
	\label{alg:hpro}
  \small
	\textbf{Input:} iterate \( \boldsymbol{X}^{1,1;\nu_1}_{l,s} \), correction \( \boldsymbol{Y}^{{h}}_{l',s-1} \) \\
	\begin{enumerate}\itemsep1pt \parskip0pt \parsep0pt
		\item[\textbf{1.}] \textbf{If} \( s \le L \), set:
		\(\displaystyle
		\boldsymbol{X}^{1,1;\nu_1+1}_{l,s} = \boldsymbol{X}^{1,1;\nu_1}_{l,s} + \boldsymbol{T}^{1,1}_{l-1,l;s-1,s} \boldsymbol{Y}^{{h}}_{l-1,s-1}.
		\)
		
		\item[\textbf{2.}] \textbf{Else}, set:
		\(\displaystyle
		\boldsymbol{X}^{1,1;\nu_1+1}_{l,s} = \boldsymbol{X}^{1,1;\nu_1}_{l,s} + \boldsymbol{T}^{1,1}_{l,l;s-1,s} \boldsymbol{Y}^{{h}}_{l,s-1}.
		\)
		
		\item[\textbf{3.}] \textbf{Return} \( \boldsymbol{X}^{1,1;\nu_1+1}_{l,s} \)
	\end{enumerate}
\end{algorithm}
\vspace*{-1ex}
\begin{algorithm}[H]
	\caption{\textsc{P-Prolongate-Correct}}
	\label{alg:ppro}
  \small
	\textbf{Input:} iterate \( \boldsymbol{X}^{k,r;\nu_1}_{l,s} \), correction \( \boldsymbol{Y}_{{p}}^{k',\frac{r}{2}} \) \\
	\begin{enumerate}\itemsep1pt \parskip0pt \parsep0pt
		\item[\textbf{1.}] \textbf{If} \( r \le K \), set:
		\(\displaystyle
		\boldsymbol{X}^{k,r;\nu_1+1}_{l,s} = \boldsymbol{X}^{k,r;\nu_1}_{l,s} + \boldsymbol{T}^{\frac{k}{2},k;\frac{r}{2},r}_{l,s} \boldsymbol{Y}_{{p}}^{\frac{k}{2},\frac{r}{2}}.
		\)
		
		\item[\textbf{2.}] \textbf{Else}, set:
		\(\displaystyle
		\boldsymbol{X}^{k,r;\nu_1+1}_{l,s} = \boldsymbol{X}^{k,r;\nu_1}_{l,s} + \boldsymbol{T}^{k,k;\frac{r}{2},r}_{l,s} \boldsymbol{Y}_{{p}}^{k,\frac{r}{2}}.
		\)
		
		\item[\textbf{3.}] \textbf{Return} \( \boldsymbol{X}^{k,r;\nu_1+1}_{l,s} \)
	\end{enumerate}
\end{algorithm}

\section{An Example of Multigrid Sequence Generation}
\label{sec:example-mgseq}
We illustrate the application of \Cref{alg:ConstructHierarchy}
and~\ref{alg:CombineHierarchies} by an example.
Table~\ref{tab:multigrid_hierarchy} summarizes the multigrid hierarchy by
listing the spatial and temporal discretization parameters for geometric and
polynomial coarsening. In space, the mesh size reduces from \(h\) to \(2h\) (Level 1) and \(4h\) (Level 0), while the polynomial degree is
reduced from 2 at Level 3 to 1 from Level 2 onward. In this example, no
geometric coarsening in time is performed and polynomial coarsening in time is
applied from Level 1 to Level 0.
\begin{table}[ht]
  \centering
  \caption{\label{tab:multigrid_hierarchy}Multigrid Hierarchy Parameters for Each Level}
  \footnotesize
  \begin{tabular}{@{}cccccl@{}}
    \toprule
    &\multicolumn{2}{c}{\textbf{Space}} &
                                          \multicolumn{2}{c}{\textbf{Time}}&\multirow{ 2}{*}{Coarsening description} \\
    \textbf{Level} & \textbf{$h$-MG} & \textbf{$p$-MG} & \textbf{$h$-MG} & \textbf{$p$-MG}& \\ \midrule
    3 & \( h \)    & \( 2 \)     & \( \tau \)    & \( 2 \)    & \multirow{ 2}{*}{Polynomial
                                                                coarsening in space} \\
    2 & \( h \)    & \( 1 \) & \( \tau \)    & \( 2 \)     & \multirow{ 2}{*}{Geometric
                                                             coarsening in space}\\
    1 & \( 2h \)   & \( 1 \) & \( \tau \)    & \( 2 \)    & \multirow{ 2}{*}{Polynomial in time,
    geometric in space}\\
    0 & \( 4h \)   & \( 1 \) & \( \tau \)    & \( 1 \) \\ \bottomrule
  \end{tabular}
\end{table}

\section{\label{Sec:CP} Convergence plots} 
\begin{figure}[htb]
  \centering
    \includegraphics{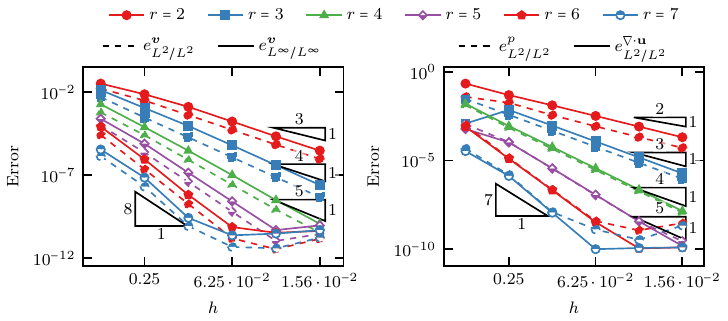}
    \caption{\label{fig:conv-stokes}Calculated errors of the velocity and pressure
      in various norms (velocity: $L^2$, $L^{\infty}$ in space-time and the
      $L^2$-norm of the divergence in space-time, pressure: $L^2$ in space-time)
      for different polynomial orders. The expected orders of convergence,
      represented by the triangles, match with the experimental orders.}
  \end{figure}

\Cref{fig:conv-stokes} shows the
convergence in various norms for all polynomial degrees and refinements
in~\Cref{sec:conv-stokes}. The iteration counts are given in Table~\ref{tab:iter-stokes}.
  \section{Pressure difference in lid-driven cavity flow}\label{sec:pdiff}
\begin{figure}\centering
    \includegraphics{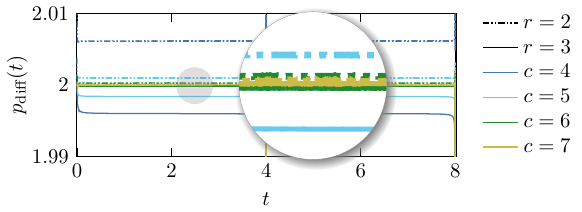}
    \caption{\label{fig:pdiff} Normalized Pressure difference \(p_{\text{diff}}(t)\) over time.}
  \end{figure}

  To assess the convergence of our discretization, we consider the normalized pressure difference 
  \[
    p_{\text{diff}}(t) = \frac{p(0.875,0.125,0.125,t) - p(0.875,0.875,0.875,t)}
    {p(0.875,0.125,0.125,t)}\,.
  \]
  We normalize the
  pressure difference in order to improve the visualization of the discretization error. In Figure~\ref{fig:pdiff} we plot \(p_{\text{diff}}(t)\) over the time interval \(I=[0,8]\).

\section*{Acknowledgments}
Computational resources (HPC cluster HSUper) have been provided by the project
hpc.bw, funded by dtec.bw - Digitalization and Technology Research Center of the
Bundeswehr. dtec.bw is funded by the European Union - NextGenerationEU.
\bibliographystyle{siamplain}
\bibliography{stokes}

\end{document}